\documentclass{amsart}
\usepackage{amsmath,amsfonts,amssymb,amscd}
\usepackage{graphics}
\usepackage{aeguill}
\usepackage{hyperref}

\newtheorem{thm}[equation]{Theorem}
\newtheorem{lem}[equation]{Lemma}
\newtheorem{prop}[equation]{Proposition}
\newtheorem{cor}[equation]{Corollary}
\theoremstyle{definition}
\newtheorem{defn}[equation]{Definition}
\newtheorem{nota}[equation]{Notation}
\newtheorem{rem}[equation]{Remark}
\newtheorem{exmp}[equation]{Example}

\DeclareMathOperator\lcm{lcm}
\DeclareMathOperator\Supp{Supp}
\DeclareMathOperator\Aut{Aut}
\newcommand{\gras}[1]{\textrm{\boldmath $#1$}} 

\newcommand{\ie}{{\em i.e.}}
\newcommand{\cf}{{\em cf.}}
\def\G{\Gamma }

\def\a{\alpha }
\def\b{\beta }
\def\l{\ell }
\def\p{\varphi }
\def\Nset{\mathbb{N}}

\title{Admissible submonoids of Artin-Tits monoids}
\author{Anatole Castella}

\begin{document}


\maketitle



\begin{abstract}We show the analogue of M\"uhlherr's [{\em Coxeter groups in Coxeter groups}, Finite Geom. and Combinatorics, Cambridge Univ. Press (1993), 277--287] for Artin-Tits monoids and for Artin-Tits groups of spherical type. That is, the submonoid (resp. subgroup) of an Artin-Tits monoid (resp. group of spherical type) induced by an {\em admissible partition\/} of the Coxeter graph is an Artin-Tits monoid (resp. group).

This generalizes and unifies the situation of the submonoid (resp. subgroup) of fixed elements of an Artin-Tits monoid (resp. group of spherical type) under the action of graph automorphisms, and the notion of {\em LCM-homomorphisms\/} defined by Crisp in [{\em Injective maps between Artin groups}, Geom. Group Theory Down Under, Canberra (1996) 119--137] and generalized by Godelle in [{\em Morphismes injectifs entre groupes d'Artin-Tits}, Algebr. Geom. Topol. \textbf{2} (2002), 519--536].

We then complete the classification of the admissible partitions for which the Coxeter graphs involved have no infinite label, started by M\"uhlherr in [{\em Some contributions to the theory of buildings based on the gate property}, Dissertation, T\"ubingen (1994)]. This leads us to the classification of Crisp's LCM-homomorphisms.

\end{abstract}



\section*{Introduction.}

In 1993-1994, M\"uhlherr introduced the notion of {\em admissible partitions\/} of a Coxeter graph to define subgroups of the associated Coxeter group that inherit a Coxeter group structure from the ambient one \cite{Mu, Mu2}. This construction generalizes the situation of the subgroup of fixed elements of a Coxeter group under the action of a group of graph automorphisms, studied by H\'ee in \cite{H}. 

The aim of this paper is to show the analogue for Artin-Tits monoids and for Artin-Tits groups of spherical type. Like in the Coxeter case, our construction generalizes the situation of the submonoid (resp. subgroup) of fixed elements of an Artin-Tits monoid (resp. group of spherical type) under the action of a group of graph automorphisms (studied in the early 2000's in \cite{Mi, C2, C2', DP}). When only finite Coxeter graphs without infinite labels are involved, our construction --- more precisely the underlying notion of morphisms between Artin-Tits monoids (or groups) --- is equivalent to the notion of {\em LCM-homomorphisms\/} defined in 1996 by Crisp \cite{C}. For arbitrary Coxeter graphs, our construction is more general than the notion of LCM-homomorphisms developed in 2002 by Godelle \cite{G}, which allowed finite Coxeter graphs with infinite labels, as it works for infinite Coxeter graphs and includes all the morphisms coming from actions of graph automorphisms and all the morphisms induced by the {\em bursts\/} of a Coxeter graph used by Paris in \cite{P}. Moreover, we show that some important combinatorial properties of those earlier defined objects (such as their respect of {\em simple elements\/} and of {\em normal forms\/}) are still valid in our more general context.

We then complete the classification of admissible partitions whose {\em type\/} has no infinite label, started by M\"uhlherr in \cite{Mu2}. With our new point of view on LCM-homomorphisms, this gives us the classification of Crisp's LCM-homomorphisms, started in \cite{C} with the notion of {\em foldings\/} of Coxeter graphs (which turn out to be nothing else but special cases of admissible partitions).

\section{Preliminaries.}

\subsection{Generalities on monoids.}\mbox{}\medskip

Let $M$ be a monoid, \ie{} a (non-empty) set endowed with an associative binary operation $M \times M \to M$, $(x,y) \mapsto xy$, with an identity element (denoted by~$1$). An element $x \in M$ is said to be a {\em left\/} (resp. {\em right}) {\em unit\/} if there exists $y \in M$ such that $xy = 1$ (resp. $yx = 1$). For example, $1$ is a left and right unit. The monoid $M$ is said to be {\em left\/} (resp. {\em right}) {\em cancellative\/} if, for all $x,\,y,\,z \in M$, $xy = xz$ (resp. $yx = zx$) implies $y = z$~; and $M$ is said to be {\em cancellative\/} if it is left and right cancellative. Note that, in a left or right cancellative monoid, left units and right units coincide.

Let $S = \{s_e \mid e \in E\}$ be a generating subset of $M$ such that the map $E \to S$, $e \mapsto s_e$, is one-to-one. A word $e_1 \cdots e_n$ on $E$ is a {\em representation\/} (on $E$) of $x \in M$ if $x = s_{e_1} \cdots s_{e_n}$, it is called {\em reduced\/} if it is of minimal length among all the representations of $x$. We denote by $\l_S(x)$ this minimal length, and call the function $\l_S : M \to \Nset$ thus defined the {\em length on $M$ with respect to $S$}.

We denote by $\preccurlyeq$ (resp. $\succcurlyeq$) the {\em left\/} (resp. {\em right}) {\em divisibility} in $M$, \ie{} for $x,\, y\in M$, we write $y \preccurlyeq x$ (resp. $x \succcurlyeq y$) if there exists $z \in M$ such that $x = yz$ (resp. $x = zy$). There are natural notions of {\em gcd}'s and {\em lcm}'s in $M$~: an element $d$ in $M$ is a {\em left gcd\/} of a non-empty subset $X \subseteq M$ if $d \preccurlyeq x$ for all $x \in X$ and if, for every $z \in M$ with this property, we get $z \preccurlyeq d$~; an element $m$ in $M$ is a {\em right lcm\/} of a non-empty subset $X \subseteq M$ if $x \preccurlyeq m$ for all $x \in M$ and if, for every $z \in M$ with this property, we get $m \preccurlyeq z$. The notions of {\em right gcd\/} and {\em left lcm\/} are defined symmetrically. If two elements $x,\, y\in M$ have a {\em unique\/} left (resp. right) lcm, we denote it by $x \vee_L y$ (resp. $x \vee_R y$)~; and if they have a {\em unique\/} left (resp. right) gcd, we denote it by $x \wedge_L y$ (resp. $x \wedge_R y$). Note that in a cancellative monoid with no non-trivial unit, gcd's and lcm's are unique when they exist.

For $x_1, \ldots ,\, x_n \in M$, we denote by $\prod_{k = 1}^nx_k$ the product $x_1 x_2 \cdots x_n$ in that order. For $x,\,y \in M$ and $n \in \Nset$, we denote by $\prod_n(x,y)$ the product $xyxy\cdots$ of $n$ terms alternatively equal to $x$ and $y$ (starting with $x$). If $M = \Nset$ endowed with the usual addition, we prefer the notation $\sum_n(x,y)$ for the sum $x+y+x+y+ \cdots$ of $n$ terms alternatively equal to $x$ and $y$ (starting with $x$).

\subsection{Generalities on Coxeter groups and Artin-Tits groups.}\label{Generalities on Coxeter and Artin-Tits groups}\mbox{}\medskip

Let $\G = (m_{i,j})_{i,j \in I}$ be a {\em Coxeter matrix\/} over an arbitrary (non necessarily finite) set $I$, \ie{} with $m_{i,j} = m_{j,i} \in \Nset_{\geqslant 1} \cup \{\infty\}$ and $m_{i,j} = 1 \Leftrightarrow i = j$. The matrix $\G$ is usually represented by its {\em Coxeter graph}, \ie{} the graph with vertex set $I$, edge set $\{\{i,j\} \mid m_{i,j} \geqslant 3\}$, and a label $m_{i,j}$ over the edge $\{i,j\}$ if $m_{i,j} \geqslant 4$. We denote by  
\[
 \begin{array}{lcl}
W_\G &=& < s_i,\, i\in I \mid s_i^2 = 1, \, \prod_{m_{i,j}}(s_i, s_j) = \prod_{m_{i,j}}(s_j, s_i), \textrm{  if  } m_{i,j}\neq \infty >, \\
B_\G &=& < \gras s_i, \, i\in I \mid \prod_{m_{i,j}}(\gras s_i,\gras s_j) = \prod_{m_{i,j}}(\gras s_j,\gras s_i), \textrm{  if  } m_{i,j}\neq \infty>, \\
B^+_\G &=& < \gras s_i, \, i\in I \mid \prod_{m_{i,j}}(\gras s_i,\gras s_j) = \prod_{m_{i,j}}(\gras s_j,\gras s_i), \textrm{  if  } m_{i,j}\neq \infty>^+,
\end{array}
\]
the {\em Coxeter group}, the {\em Artin-Tits group\/} and the {\em Artin-Tits monoid\/} associated with $\G$ respectively. Note that we may use the same symbols for the generators of $B_\G$ and $B^+_\G$ since Paris showed in \cite{P} that $B^+_\G$ identifies with the submonoid of $B_\G$ generated by the $\gras s_i$, $i \in I$ (he actually proved this result when $I$ is finite, but this implies the general case). Set $S_\G = \{s_i\mid i\in I\}$ and $\gras S_\G = \{\gras s_i\mid i\in I\}$~; we say that the pair $(W_\G,S_\G)$ (resp. $(B_\G,\gras S_\G)$, resp. $(B^+_\G,\gras S_\G)$) is the {\em Coxeter\/} (resp. {\em Artin-Tits}, resp. {\em positive Artin-Tits}) {\em system\/} of {\em type\/} $\G$. Note that $W_\G$ is generated by $S_\G$ as a monoid. We denote by the same letter $\l$ the lengths on $W_\G$ with respect to $S_\G$, and on $B^+_\G$ with respect to $\gras S_\G$, and call them {\em standard\/} lengths.

Let $\G = (m_{i,j})_{i,j\in I}$ and $\G' = (m'_{i',j'})_{i',j'\in I'}$ be two Coxeter matrices. An isomorphism from $\G$ onto $\G'$ is a bijective map $f : I  \to I'$ such that $m_{i,j} = m'_{f(i),f(j)}$ for all $i,\, j \in I$. In particular, we denote by $\Aut(\G)$ the automorphism group of $\G$. We say that two pairs $(G_1,S_1)$ and $(G_2,S_2)$, where $G_i$ is a group (resp. a monoid) generated by $S_i$ ($i = 1,\, 2$), are {\em isomorphic\/} if there exists an isomorphism $f : G_1 \to G_2$ that maps $S_1$ onto $S_2$. For example, the two systems $(W_\G,S_\G)$ and $(W_{\G'},S_{\G'})$ (resp. $(B_\G,\gras S_\G)$ and $(B_{\G'},\gras S_{\G'})$, resp. $(B^+_\G,\gras S_\G)$ and $(B^+_{\G'},\gras S_{\G'})$) are isomorphic if and only if so are $\G$ and $\G'$. 

\subsubsection{Simple elements.}\mbox{}\medskip

Let $\pi_\G : B_\G \to W_\G$ be the canonical morphism sending $\gras s_i$ on $s_i$ for all $i \in I$.

The order of $s_is_j$ in $W_\G$ is exactly $m_{i,j}$ \cite[Ch.~V, n$^\circ$~4.3, Prop.~4]{B}. In particular, the map $I \to S_\G$, $i \mapsto s_i$, and hence the map $I \to \gras S_\G$, $i \mapsto \gras s_i$, are one-to-one. Tits showed in \cite[Thm.~3]{T} that two reduced representations on $I$ of an element $w \in W$ only differ from a finite sequence of transformations --- called {\em braid relations} --- of the form $\prod_{m_{i,j}}(i,j) \leadsto \prod_{m_{i,j}}(j,i)$ with $i,\, j \in I$ such that $i\neq j$ and $m_{i,j} \neq \infty$. This property makes the following definition allowable~:

\begin{defn}[simple elements]The canonical morphism $\pi_\G : B_\G \to W_\G$ has a section $w \mapsto \gras w \in B^+_\G$ where $\gras w$ is represented on $I$ by one (and hence any) reduced representation of $w$ on $I$. We say that such an element $\gras w$ in $B^+_\G$ is {\em simple\/} and set $\gras W_\G = \{\gras w \mid w\in W_\G\} = \{x\in B^+_\G \mid \l(x) = \l(\pi_\G(x))\}$.
\end{defn}

\subsubsection{Standard parabolicity, sphericity and irreducibility.}\mbox{}\medskip

Let $J \subseteq I$. We set $\G_J = (m_{i,j})_{i,j\in J}$ (it is a Coxeter matrix)~; and we denote by $W_J$ (resp. $B_J$, resp. $B^+_J$) the subgroup of $W_\G$ (resp. the subgroup of $B_\G$, resp. the submonoid of $B^+_\G$) generated by $\{s_j \mid j \in J\}$ (resp. $\{\gras s_j \mid j \in J\}$). 

\begin{defn}[standard parabolicity]The subgroups $W_J$ (resp. subgroups $B_J$, resp. submonoids $B^+_J$), $J \subseteq I$, of $W_\G$ (resp. $B_\G$, resp. $B^+_\G$) are called {\em standard parabolic} ({\em with respect to $\G$}).
\end{defn}

The pair $(W_J,\{s_j \mid j\in J\})$ (resp.  $(B_J,\{\gras s_j \mid j\in J\})$, resp. $(B^+_J,\{\gras s_j \mid j\in J\})$) is (isomorphic to) the Coxeter (resp. Artin-Tits, resp. positive Artin-Tits) system of type $\G_J$ (see \cite[Ch.~IV, n$^\circ$~1.8, Thm.~2]{B} for the Coxeter case, \cite[Ch.~II, Thm.~4.13]{vdL} for the Artin-Tits case with $I$ finite --- which implies the general result ---, the positive Artin-Tits case being obvious).

Moreover, the standard length on $W_J$ (resp. $B^+_J$) is induced by the one on $W_\G$ (resp. $B^+_\G$) \cite[Ch.~IV, n$^\circ$~1.8, Cor.~4]{B}. This implies that $\gras W_J = \gras W_\G \cap B_J^+$.

\begin{defn}[sphericity]The Coxeter matrix $\G_J$ is called {\em spherical} --- and the subset $J$ of $I$ is called {\em spherical} ({\em with respect to $\G$}) --- if $W_J$ is finite.  In that case, the subgroups $W_J$, $B_J$, and submonoid $B^+_J$ are also called {\em spherical}. 
\end{defn}

In a finite Coxeter group, there exists a unique element of maximal standard length, which is of order two if not trivial \cite[Ch.~IV, \S~1, Ex.~22]{B}. If $J$ is spherical, we denote by $r_J$ the unique element of maximal standard length in $W_J$ and by $\gras r_J$ its image in $\gras W_J$ (\ie{} the unique element of maximal standard length in $\gras W_J$).

\begin{defn}[irreducibility]The matrix $\G$ is said to be {\em reducible\/} if there exists a partition of cardinality two $\{J,K\}$ of $I$ such that $m_{j,k} = 2$ for every pair $(j,k) \in J \times K$. In that case, we write $\G = \G_J \times \G_K$, as we have $W_\G = W_J \times W_K$, $B_\G = B_J \times B_K$ and $B^+_\G = B^+_J \times B^+_K$. If this is not the case, then $\G$ is said to be {\em irreducible}~; this is precisely when the Coxeter graph of $\G$ is connected.
\end{defn}

We assume that the reader is familiar with the list of the irreducible spherical Coxeter graphs, which can be found for example in \cite[Ch.~VI, n$^\circ$~4.1, Thm.~1]{B}.

\subsubsection{Properties of $B^+_\G$.}\mbox{}\medskip

Since the defining relations of $B^+_\G$ are homogeneous, the standard length of $B^+_\G$ is {\em additive}, \ie{} $\l(xy) = \l(x) + \l(y)$ for all $x,\,y \in B^+_\G$. This clearly implies that $B^+_\G$ has no non-trivial unit. Moreover, $B^+_\G$  is cancellative \cite[Prop.~2.3]{BS} (hence gcd's and lcm's are unique when they exist), and two elements of $B^+_\G$ always have left and right gcd's, and have a right (resp. left) lcm as soon as they have a right (resp. left) common multiple \cite[Props.~4.1 and 4.2]{BS}.

\begin{exmp}Let $J \subseteq I$ be non-empty. By {\cite[Thm.~5.6]{BS}}, the elements $\gras s_j$, $j\in J$, have a (left or right) lcm if and only if $\G_J$ is spherical, and in that case their (left and right) lcm is $\gras r_J$ {\cite[Prop.~5.7]{BS}}. In particular, two elements $\gras s_i$ and $\gras s_j$ have a (left or right) lcm if and only if $m_{i,j} \neq \infty$, in which case $\gras s_i \vee_R \gras s_j = \gras s_i \vee_L \gras s_j =\gras r_{\{i,j\}} = \prod_{m_{i,j}}(\gras s_i,\gras s_j) = \prod_{m_{i,j}}(\gras s_j,\gras s_i)$.
\end{exmp}

In \cite[Prop.~2.1]{Mi}, Michel showed that for all $x \in B^+_\G$, there exists a unique maximal (for $\preccurlyeq$) element $L(x)$ in the set $\{\gras w \in \gras W_\G \mid \gras w \preccurlyeq x\}$ of all simple left divisors of $x$. The maximal simple right divisor $R(x)$ of $x$ is defined symmetrically.

\begin{defn}[normal forms]\label{defs formes normales}
The {\em left normal form\/} of a non-trivial element $x \in B^+_\G$ is the unique sequence $(x_1,\ldots,\,x_n)$ of elements of $\gras W_\G$ such that $x = x_1\cdots x_n$, $x_n \neq 1$ and $x_k = L(x_kx_{k+1} \cdots x_n)$ for $1\leqslant k\leqslant n-1$. {\em Right normal forms\/} are defined symmetrically.
\end{defn}

It is clear that $B^+_\G$ generates $B_\G$ (as a group). If $\G$ is spherical, $B_\G$ is more precisely the {\em group of fractions\/} of $B^+_\G$, \ie{} every $b \in B_\G$ can be written $b = x^{-1}y = x'y'^{-1}$ for $x,\, y,\, x',\, y' \in B^+_\G$ \cite[Prop.~5.5]{BS}. 

\begin{defn}[irreducible fractions]\label{def irred left form}Assume that $\G$ is spherical and fix $b \in B_\G$. Then { \cite[Cor.~7.5]{DP}\/} shows that there exists a {\em unique\/} pair $(x,y)$ (resp. $(x',y')$) in $(B^+_\G)^2$ such that $b = x^{-1}y$ and $x\wedge_L y =1$ (resp. $b = x'y'^{-1}$ and $x'\wedge_R y' =1$). We say that this pair $(x,y)$ (resp. $(x',y')$) is an {\em irreducible left\/} (resp. {\em right}) {\em fraction}, and is the {\em irreducible left\/} (resp. {\em right}) {\em form\/} of $b$.
\end{defn}

\section{Admissible partitions --- The work of M\"uhlherr.}

In this section, we recall the definition of an admissible partition of a Coxeter graph and the principal results of \cite{Mu} on the subgroup of the associated Coxeter group defined by such a partition. Let  $\G = (m_{i,j})_{i,j\in I}$ be a Coxeter matrix and let $W = W_\G$.

\subsection{Definitions.}\label{definitions des partitions admissibles}\mbox{}\medskip

\begin{defn}[{\cite{Mu}}]\label{defs sph�riques etc}We say that a partition $\tilde I$ of $I$ is {\em spherical} ({\em with respect to $\G$}) --- or by abuse of language is a {\em spherical partition of $\G$} --- if, for all $\a \in \tilde I$, $\G_\a$ is spherical (\ie{} $W_\a$ is finite). In that case, we denote by 
\begin{itemize}
\item $\tilde S = \{r_\a \mid \a \in \tilde I\}$ the set of all $r_\a$, $\a \in \tilde I$ (recall that $r_\a$ is the unique element of maximal standard length in $W_\a$),
\item $\tilde W = \langle \tilde S \rangle$ the subgroup of $W$ generated by $\tilde S$,
\item $\tilde l = l_{\tilde S}$ the length on $\tilde W$ with respect to $\tilde S$ ($\tilde W$ is generated by $\tilde S$ as a monoid),
\item $\tilde \G = ( | r_\a r_\b|)_{\a,\b \in \tilde I}$ the Coxeter matrix of orders of the products $r_\a r_\b$ in $W$. We call $\tilde \G$ the {\em type} of $\tilde I$. 
\end{itemize}
Moreover for $\a_1,\ldots,\, \a_n \in \tilde I$ and $w = \prod_{k=1}^nr_{a_k} \in \tilde W$, we say that the word $\prod_{k=1}^n\a_k$ on $\tilde I$ is {\em compatible} --- or is a {\em compatible representation of\/} $w$ --- ({\em with respect to $\G$}), if $\l(w) = \sum_{k = 1}^n\l(r_{\a_k})$. 
\end{defn}

Note that we always have $\l(w) \leqslant \sum_{k = 1}^n\l(r_{\a_k})$, and the equality holds precisely when the representation $R_{\a_1}\cdots R_{\a_n}$ of $w$ on $I$, where for $1 \leqslant k \leqslant n$ the word $R_{\a_k}$ is a reduced representation of $r_{\a_k}$ on $I$, is reduced.

\begin{nota}Let $w \in W$. We set $\begin{cases} I^+(w) = \{i\in I \mid \l(ws_i) = \l(w)+1\}\\ I^-(w) = \{i\in I \mid \l(ws_i) = \l(w)-1\} \end{cases}$.
\end{nota}

Note that $I^-(w)$ is a spherical subset of $I$ \cite[Lem.~2.8]{Mu}.

\begin{defn}[{\cite{Mu}}]\label{def admissibilit�}Let $\tilde I$ be a partition of $I$. We say that $\tilde I$ is {\em admissible} ({\em with respect to $\G$}) --- or by abuse of language is an {\em admissible partition of $\G$} --- if it is a spherical partition of $\G$ such that, for all $(w,\a) \in \tilde W \times \tilde I$, either $\a \subseteq I^+(w)$ or $\a \subseteq I^-(w)$.
\end{defn}

\begin{rem}Let $\a$ be a spherical subset of $I$ and $w \in W$. Then $\a \subseteq I^-(w)$ (resp. $\a \subseteq I^+(w)$) if and only if $\l(wr_\a) = \l(w)-\l(r_\a)$ (resp. $\l(wr_\a) = \l(w)+\l(r_\a)$) { \cite[Lems.~2.4 and 2.8]{Mu}}.
\end{rem}

\subsection{Admissible partitions and Coxeter groups.}\label{Partitions admissibles et groupes de Coxeter.}\mbox{}\medskip

The two main results of \cite{Mu} are the following theorems~: 

\begin{thm}[{\cite[Thm.~1.1]{Mu}}]\label{add implique inj}Let $\tilde I$ be an admissible partition of $\G$, of type $\tilde \G$. Then the pair $(\tilde W, \tilde S)$ is (isomorphic to) the Coxeter system of type $\tilde \G$.
\end{thm}

\begin{thm}[{\cite[Thm.~1.2]{Mu}}]\label{admissibilit� locale}Let $\tilde I$ be a partition of $\G$. The following conditions are equivalent~: 
\begin{enumerate}
\item $\tilde I$ is an admissible partition of $\G$,
\item for all $\a,\,\b \in \tilde I$ with $\a \neq \b$, $\{\a,\b\}$ is an admissible partition of $\G_{\a \cup \b}$.
\end{enumerate}
\end{thm}

So proving the admissibility of a partition reduces to proving the admissibility of partitions of cardinality two. The following lemma gives a criterion for that. It is left as an exercise in \cite{Mu}, but for convenience and because it will be of great importance for our purpose, we prove it below, following \cite{D}. Note that our condition {(\ref{caract ad3})\/} is slightly weaker than the one of \cite[Lem.~3.3]{Mu}~; this formulation simplifies the proof of the second part of the lemma and will be useful later in section \ref{LCM-homoms}. From now on, we call {\em $2$-partition\/} a partition of cardinality two. 

\begin{lem}[{\cite[Lem.~3.3]{Mu}}]\label{caract ad}Let $\tilde I = \{\a,\b\}$ be a spherical $2$-partition of $\G$. 
\begin{enumerate}
 \item The following conditions are equivalent~: 
\begin{enumerate}
\item \label{caract ad1} $\tilde I$ is an admissible partition of $\G$,
\item \label{caract ad3} for every integer $0 \leqslant n < |r_\a r_\b| + 1$, the words $\prod_n(\a, \b)$ and $\prod_n(\b, \a)$ are compatible.
\end{enumerate}
\item If $\G$ is spherical, then $|r_\a r_\b| \neq \infty$ and conditions {\em (\ref{caract ad1})} and {\em (\ref{caract ad3})} above are equivalent to the following condition~:
\begin{enumerate}
\item \label{caract ad cas spherique2} the words $\prod_{|r_\a r_\b|}(\a, \b)$ and $\prod_{|r_\a r_\b|}(\b, \a)$ are compatible.
\end{enumerate}
Moreover, we get in that case $\prod_{|r_\a r_\b|}(r_\a, r_\b) = \prod_{|r_\a r_\b|}(r_\b, r_\a) = r_I$.
\end{enumerate}

\end{lem}

\proof  The subgroup $\tilde W = \langle r_\a,r_\b\rangle$ of $W$ is a dihedral group of order $2|r_\a r_\b|$, hence the reduced representations on $\tilde I$ of the elements of $\tilde W$ are the words $\prod_n(\a, \b)$ and $\prod_n(\b, \a)$ for every integer $0 \leqslant n < |r_\a r_\b|+1$. 

Suppose {(\ref{caract ad3})\/} and let us show {(\ref{caract ad1})}. Let $w = \prod_n(r_\a,r_\b) \in \tilde W$ for some $0 \leqslant n < |r_\a r_\b|+1$. We have to show that either $\a \subseteq I^+(w)$ or $\a \subseteq I^-(w)$, and the same for $\b$. We can assume $w \neq 1$ (because $\a \cup \b = I = I^+(1)$). For $k \in \Nset$, set $\a_k = \a$ if $k$ is odd and $\a_k = \b$ if $k$ is even. Since $\prod_n(\a,\b)$ is compatible, we get $\a_n \subseteq I^-(w)$. If $|r_\a r_\b| \neq \infty$ and if $n = |r_\a r_\b|$, we thus get by symmetry $\a \cup \b = I = I^-(w)$.  If $n < |r_\a r_\b|$, then the word $\prod_{n+1}(\a,\b)$ is compatible, whence $\a_{n+1} \subseteq I^-(wr_{\a_{n+1}})$ and hence $\a_{n+1} \subseteq I^+(w)$.

Suppose {(\ref{caract ad1})\/} and let us show {(\ref{caract ad3})}. We first prove, by induction on $\l(w)$, that every $w \in \tilde W$ admits a compatible representation on $\tilde I$. If $w = 1$ this is obvious, else let $i \in I$ be such that $\l (ws_i)  = \l(w) -1$. There is no loss of generality in assuming that $i \in \a$. Since $\tilde I$ is admissible, we have $\a \subseteq I^-(w)$, and $\l(wr_\a) = \l(w)-\l(r_\a)$. By induction, $wr_\a$ admits a compatible representation $\a_1 \cdots \a_n$, and $\a_1 \cdots \a_n\a$ is then a compatible representation of $w$. 

Now, fix an integer $0 \leqslant n < |r_\a r_\b|+1$ and consider the word $\prod_{n}(\a,\b)$. If $n < |r_\a r_\b|$, then this word is the unique reduced representation on $\tilde I$ of the element $w = \prod_n(r_\a,r_\b) \in \tilde W$, so it must be the existing compatible representation of $w$ (it is clear that a non-reduced word on $\tilde I$ cannot be compatible). It remains to prove that, if $|r_\a r_\b| \neq \infty$ and if $\prod_{|r_\a r_\b|}(\a,\b)$ is compatible, then so is $\prod_{|r_\a r_\b|}(\b,\a)$. This is clear if $|r_\a r_\b|$ is even, so assume $|r_\a r_\b|$ odd and set $w = \prod_{|r_\a r_\b|}(r_\a,r_\b) = \prod_{|r_\a r_\b|}(r_\b,r_\a)$ and $w' = \prod_{|r_\a r_\b|-1}(r_\b,r_\a)$. The word $\prod_{|r_\a r_\b|-1}(\b,\a)$ is the unique reduced representation of $w'$, hence it is compatible and we have $\a \subseteq I^-(w')$. Since $w'$ is not the element of maximal standard length in $W$, we get $\b \not \subseteq I^-(w')$, whence $\b \subseteq I^+(w')$ by admissibility, and hence $\prod_{|r_\a r_\b|}(\b,\a)$ is a compatible representation of $w$. 

If $\G$ is spherical, then it is clear that $|r_\a r_\b| \neq \infty$ and {(\ref{caract ad3})\/} implies {(\ref{caract ad cas spherique2})}. Conversely, if {(\ref{caract ad cas spherique2})\/} holds, then  for all $0 \leqslant n < |r_\a r_\b|$, the prefix $\prod_n(\a, \b)$ of $\prod_{|r_\a r_\b|}(\a, \b)$ (resp. $\prod_n(\b, \a)$ of $\prod_{|r_\a r_\b|}(\b, \a)$) is necessarily compatible, whence {(\ref{caract ad3})}. Now consider $w = \prod_{|r_\a r_\b|}(r_\a, r_\b) = \prod_{|r_\a r_\b|}(r_\b, r_\a)$ in $\tilde W$. Since both words $\prod_{|r_\a r_\b|}(\a, \b)$ and $\prod_{|r_\a r_\b|}(\b, \a)$ are compatible, we get $\a \cup \b = I = I^-(w)$, whence $w = r_I$. 
\qed\medskip

Let us conclude this subsection with some further properties of admissible partitions~:

\begin{prop}[{\cite[Prop.~3.5, A1]{Mu} and \cite[Lem.~2.5.5]{Mu2}}]\label{I spherique ssi tildeI sph�rique}Let $\tilde I$  be an admissible partition of $\G$, of type $\tilde \G$, and let $w \in \tilde W$. 
\begin{enumerate}
\item a representation of $w$ on $\tilde I$ is reduced if and only if it is compatible,
\item $\G$ is spherical if and only if so is $\tilde \G$, in which case $r_I = r_{\tilde I}$. 
\end{enumerate}
\end{prop}

\subsection{Examples.}\label{exemples de partitions admissibles}\mbox{}\medskip

Let $\G = (m_{i,j})_{i,j \in I}$ be a Coxeter matrix and $G$ be a subgroup of $\Aut(\G)$. The action of $G$ on $I$ induces an action of $G$ on $W_\G$ which preserves the standard length. If $\a$ is an orbit of $I$ under $G$, then $G$ stabilizes $W_\a$ and hence, if $\a$ is spherical, $G$ fixes $r_\a$ (which is the unique element of maximal standard length in $W_\a$). So if we denote by $\tilde J$ the set of spherical orbits of $I$ under $G$, by $J= \bigcup_{\a \in \tilde J} \a \subseteq I$ their union and if we set $\tilde S = \{ r_\a \mid \a \in \tilde J \}$ and $\tilde W = \langle \tilde S \rangle$, we get that $\tilde W$ is included in the subgroup $(W_\G)^G$ of fixed points of $W_\G$ under $G$, and that $\tilde J$ is an admissible partition of $\G_J$. Let $\tilde \G$ be the type of $\tilde J$.

In fact, it can be shown that $\tilde W = (W_\G)^G$, hence $((W_\G)^G, \tilde S)$ is (isomorphic to) the Coxeter system of type $\tilde \G$ \cite[Thm.~1.3]{Mu}. See \cite[Cor.~3.5]{H} for the original proof of that result.

\begin{exmp}\label{sym spheriques}Here are symbolized the non-trivial automorphisms of the spherical irreducible Coxeter graphs, and the type of the different sets of orbits we get (see { \cite[section~2.5]{Mu2}} or section { \ref{Classification}} below for justifications)~: 
\begin{center}\setlength{\unitlength}{1.3pt}
\begin{picture}(85,50)(-25,-30)
\put(-25,5){ $A_{2n-\!1}$}
\put(-25,-15){\footnotesize  $(n\!\geqslant \!2)$}

\put(0,0){\circle{4}}
\put(15,0){\circle*{4}}
\put(30,0){\circle{4}}\put(30,0){\circle{6}}
\put(45,7.5){\circle*{4}}\put(45,7.5){\circle{6}}
\put(0,15){\circle{4}}
\put(15,15){\circle*{4}}
\put(30,15){\circle{4}}\put(30,15){\circle{6}}

\put(2,0){\line(1,0){11}}
\put(2,15){\line(1,0){11}}
\put(17.5,0){\small{$\ldots$}}
\put(17.5,15){$\ldots$}
\put(32,15){\line(2,-1){14}}
\put(32,0){\line(2,1){14}}

\put(22.5,-8){\vector(0,-1){10}}

\put(-25,-28){$B_n$}
\put(0,-25){\circle{4}}
\put(15,-25){\circle*{4}}
\put(30,-25){\circle{4}}\put(30,-25){\circle{6}}
\put(45,-25){\circle*{4}}\put(45,-25){\circle{6}}

\put(2,-25){\line(1,0){11}}
\put(17.5,-25){$\ldots$}
\put(32,-25){\line(1,0){11}}\put(35,-23){\small{$4$}}
\end{picture}
\begin{picture}(85,50)(-25,-30)
\put(-25,5){$A_{2n}$}
\put(-25,-15){\footnotesize $(n\!\geqslant \! 2)$}

\put(0,0){\circle{4}}
\put(15,0){\circle*{4}}
\put(30,0){\circle{4}}\put(30,0){\circle{6}}
\put(45,0){\circle*{4}}\put(45,0){\circle{6}}
\put(0,15){\circle{4}}
\put(15,15){\circle*{4}}
\put(30,15){\circle{4}}\put(30,15){\circle{6}}
\put(45,15){\circle*{4}}\put(45,15){\circle{6}}

\put(2,0){\line(1,0){11}}
\put(2,15){\line(1,0){11}}
\put(17.5,0){$\ldots$}
\put(17.5,15){$\ldots$}
\put(32,0){\line(1,0){11}}
\put(32,15){\line(1,0){11}}
\put(45,0){\line(0,1){14}}

\put(22.5,-8){\vector(0,-1){10}}

\put(-25,-28){$B_n$}
\put(0,-25){\circle{4}}
\put(15,-25){\circle*{4}}
\put(30,-25){\circle{4}}\put(30,-25){\circle{6}}
\put(45,-25){\circle*{4}}\put(45,-25){\circle{6}}

\put(2,-25){\line(1,0){11}}
\put(17.5,-25){$\ldots$}
\put(32,-25){\line(1,0){11}}\put(35,-23){\small{$4$}}
\end{picture}
\begin{picture}(85,50)(-25,-30)
\put(-25,5){$D_{n+1}$}
\put(-25,-15){\footnotesize  $(n\!\geqslant \! 3)$}

\put(0,7.5){\circle{4}}
\put(15,7.5){\circle*{4}}
\put(30,7.5){\circle{4}}\put(30,7.5){\circle{6}}
\put(45,0){\circle*{4}}\put(45,0){\circle{6}}
\put(45,15){\circle*{4}}\put(45,15){\circle{6}}

\put(2,7.5){\line(1,0){11}}
\put(17.5,7.5){$\ldots$}
\put(31.5,9){\line(2,1){14}}
\put(31.5,6){\line(2,-1){14}}

\put(22.5,-8){\vector(0,-1){10}}

\put(-25,-28){$B_n$}
\put(0,-25){\circle{4}}
\put(15,-25){\circle*{4}}
\put(30,-25){\circle{4}}\put(30,-25){\circle{6}}
\put(45,-25){\circle*{4}}\put(45,-25){\circle{6}}

\put(2,-25){\line(1,0){11}}
\put(17.5,-25){$\ldots$}
\put(32,-25){\line(1,0){11}}\put(35,-23){\small{$4$}}
\end{picture}
\end{center}

\begin{center} --- \end{center}

\begin{center}\setlength{\unitlength}{1.3pt}
\begin{picture}(65,50)(-25,-30)
\put(-25,5){$D_4$}

\put(0,7.5){\circle{4}}
\put(15,7.5){\circle*{4}}
\put(15,0){\circle*{4}}
\put(15,15){\circle*{4}}

\put(2,7.5){\line(1,0){14}}
\put(1.5,9){\line(2,1){14}}
\put(1.5,6){\line(2,-1){14}}

\put(7.5,-5){\vector(0,-1){10}}

\put(-25,-28){$G_2$}
\put(0,-25){\circle{4}}
\put(15,-25){\circle*{4}}
\put(2,-25){\line(1,0){11}}\put(5,-23){\small{$6$}}
\end{picture}
\begin{picture}(90,50)(-25,-30)
\put(-25,5){$E_6$}

\put(0,7.5){\circle{4}}
\put(15,7.5){\circle*{4}}
\put(30,0){\circle{4}}\put(30,0){\circle{6}}
\put(30,15){\circle{4}}\put(30,15){\circle{6}}
\put(45,0){\circle*{4}}\put(45,0){\circle{6}}
\put(45,15){\circle*{4}}\put(45,15){\circle{6}}

\put(2,7.5){\line(1,0){11}}
\put(16.5,9){\line(2,1){11.5}}
\put(16.5,6){\line(2,-1){11.5}}
\put(32,0){\line(1,0){14}}
\put(32,15){\line(1,0){14}}

\put(22.5,-5){\vector(0,-1){10}}

\put(-25,-28){$F_4$}
\put(0,-25){\circle{4}}
\put(15,-25){\circle*{4}}
\put(30,-25){\circle{4}}\put(30,-25){\circle{6}}
\put(45,-25){\circle*{4}}\put(45,-25){\circle{6}}

\put(2,-25){\line(1,0){11}}
\put(17,-25){\line(1,0){11}}\put(20,-23){\small{$4$}}
\put(32,-25){\line(1,0){11}}
\end{picture}
		\begin{picture}(65,50)(-25,-30)
\put(-25,5){$F_4$}
\put(0,0){\circle{4}}
\put(0,15){\circle{4}}
\put(15,0){\circle*{4}}
\put(15,15){\circle*{4}}

\put(2,0){\line(1,0){12}}
\put(2,15){\line(1,0){12}}
\put(15,0){\line(0,1){14}}\put(17,5){\small{$4$}}

\put(7.5,-5){\vector(0,-1){10}}

\put(-25,-28){$I_2$\small{$(8)$}}
\put(0,-25){\circle{4}}
\put(15,-25){\circle*{4}}
\put(2,-25){\line(1,0){11}}\put(5,-23){\small{$8$}}
\end{picture}
\begin{picture}(30,50)(-25,-30)
		\put(-25,5){$I_2$\small{$(m)$}}
		\put(-25,-15){\footnotesize $(m\!\geqslant \!3)$}
\put(5,0){\circle*{4}}
\put(5,15){\circle*{4}}
\put(5,0){\line(0,1){14}}\put(7,5){\small{$m$}}

\put(5,-8){\vector(0,-1){10}}

\put(-25,-28){$A_1$}
\put(5,-25){\circle*{4}}
\end{picture}
\end{center}
\end{exmp}

\begin{exmp}\label{H3 dans D6}Here are two admissible partitions that are not the set of orbits of an action of graph automorphisms (see { \cite[section~2.5]{Mu2}}, subsection { \ref{The bursts of a Coxeter graph.}\/} or section { \ref{Classification}\/} below for justifications)~:
\begin{center}\setlength{\unitlength}{1.3pt}
\begin{picture}(95,50)(-25,-30)
\put(-25,5){$D_6$}

\put(0,0){\circle{4}}
\put(0,15){\circle{4}}
\put(15,0){\circle*{4}}
\put(15,15){\circle*{4}}
\put(30,0){\circle{4}}\put(30,0){\circle{6}}
\put(30,15){\circle{4}}\put(30,15){\circle{6}}

\put(2,0){\line(1,0){11.1}}
\put(2,15){\line(1,0){11.1}}
\put(1.5,1.5){\line(1,1){13}}
\put(17,0){\line(1,0){11}}
\put(17,15){\line(1,0){11}}

\put(15,-8){\vector(0,-1){10}}

\put(-25,-28){$H_3$}
\put(0,-25){\circle{4}}
\put(15,-25){\circle*{4}}
\put(30,-25){\circle{4}}\put(30,-25){\circle{6}}

\put(2,-25){\line(1,0){11.1}}\put(5,-23){\small{$5$}}
\put(17,-25){\line(1,0){11}}
\end{picture}
\begin{picture}(95,50)(-25,-30)
\put(-25,5){$E_8$}

\put(0,0){\circle{4}}
\put(0,15){\circle{4}}
\put(15,0){\circle*{4}}
\put(15,15){\circle*{4}}
\put(30,0){\circle{4}}\put(30,0){\circle{6}}
\put(30,15){\circle{4}}\put(30,15){\circle{6}}
\put(45,0){\circle*{4}}\put(45,0){\circle{6}}
\put(45,15){\circle*{4}}\put(45,15){\circle{6}}

\put(2,0){\line(1,0){11.1}}
\put(2,15){\line(1,0){11.1}}
\put(1.5,1.5){\line(1,1){13}}
\put(17,0){\line(1,0){11}}
\put(17,15){\line(1,0){11}}
\put(32,0){\line(1,0){11}}
\put(32,15){\line(1,0){11}}

\put(22.5,-8){\vector(0,-1){10}}

\put(-25,-28){$H_4$}
\put(0,-25){\circle{4}}
\put(15,-25){\circle*{4}}
\put(30,-25){\circle{4}}\put(30,-25){\circle{6}}
\put(45,-25){\circle*{4}}\put(45,-25){\circle{6}}

\put(2,-25){\line(1,0){11.1}}\put(5,-23){\small{$5$}}
\put(17,-25){\line(1,0){11}}
\put(32,-25){\line(1,0){11}}
\end{picture}
\end{center}
\end{exmp}

\section{Admissible partitions and Artin-Tits monoids or groups.}\label{LCM-homoms}

In subsection \ref{Morphismes admissibles entre mono�des d'Artin-Tits.} below, we introduce the submonoid of an Artin-Tits monoid (resp. the subgroup of an Artin-Tits group), and the morphism between Artin-Tits monoids or groups, induced by an admissible partition of a Coxeter graph, and we establish the analogue of \cite[Thm.~1.1]{Mu} (\cf{} theorem \ref{add implique inj} above) for Artin-Tits monoids and for Artin-Tits groups of spherical type. 

In subsection \ref{Les morphismes admissibles dans la litt�rature.}, we explain how our constructions generalize the situations of the submonoids (resp. subgroups) of fixed elements of an Artin-Tits monoid (resp. group of spherical type) under the action of graph automorphisms, of the LCM-homomorphisms \cite{C,G}, and of the morphisms between Artin-Tits monoids (or groups) induced by the {\em bursts\/} of a Coxeter graph \cite{P}.

In subsection \ref{Quelques propri�t�s des morphismes admissibles.}, we show that some important properties of submonoids of fixed elements of an Artin-Tits monoid under the action of graph automorphisms and of LCM-homomorphisms extend to our settings. In particular, we establish them for the morphisms induced by the {\em bursts\/} of a Coxeter graph \cite{P}, for which they were not known when Coxeter graphs with infinite labels are involved.

But let us begin this section by recalling the notion of {\em morphisms that respect lcm's\/} defined by Crisp in \cite{C}. It is the key-tool in the proofs of the injectivity of the LCM-homomorphisms in \cite{C,G}, and plays a similar role for our main result of subsection \ref{Morphismes admissibles entre mono�des d'Artin-Tits.}.

\subsection{Morphisms that respect lcm's.}\label{Morphismes qui respectent les ppcm.}\mbox{}\medskip

Let $\G = (m_{i,j})_{i,j\in I}$ and $\tilde \G = (\tilde m_{\a,\b})_{\a,\b\in \tilde I}$ be two Coxeter matrices (where $\tilde I$ is here an arbitrary set). If $x$ and $y$ are two elements of $B^+_\G$ (resp. $B^+_{\tilde \G}$), we say for short that $x \vee_R y$ {\em exists\/} in $B^+_\G$ (resp. $B^+_{\tilde \G}$) to state that $x$ and $y$ admit a right lcm in $B^+_\G$  (resp. $B^+_{\tilde \G}$). 

\begin{defn}[{\cite[Def.~1.1]{C}}]We say that a morphism $\p : B^+_{\tilde \G} \to B^+_\G$ {\em respects right lcm's\/} if~: 
\begin{enumerate}
\item for all $\a \in \tilde I$, $\p(\gras s_{\a}) \neq 1$,
\item for all $\a,\, \b \in \tilde I$, $\gras s_{\a} \vee_R \gras s_{\b}$ exists in $B^+_{\tilde \G}$ if and only if $\p(\gras s_{\a}) \vee_R \p(\gras s_{\b})$ exists in $B^+_\G$, in which case $\p(\gras s_{\a}) \vee_R \p(\gras s_{\b}) = \p(\gras s_{\a} \vee_R \gras s_{\b})$.
\end{enumerate}

Morphisms that {\em respect left lcm's\/} are defined symmetrically, and we say that such a morphism {\em respects lcm's\/} if it respects right and left lcm's.
\end{defn}

\begin{prop}[{\cite[Thm.~8]{C2}}]\label{props morphisme ppcm}Let $\p : B^+_{\tilde \G} \to B^+_\G$ a morphism that respects right lcm's. Then~:
\begin{enumerate}
\item for all $x,\,y \in B^+_{\tilde \G}$, $x \vee_R y$ exists in $B^+_{\tilde \G}$ if and only if $\p(x) \vee_R \p(y)$ exists in $B^+_{\G}$, in which case $\p(x) \vee_R \p(y) = \p(x\vee_R y)$,
\item for all $x,\,y \in B^+$, $\p(x) \preccurlyeq \p(y) \Rightarrow x \preccurlyeq y$. In particular, $\p$ is injective. 
\end{enumerate}
\end{prop}

Of course, the symmetrical version of proposition \ref{props morphisme ppcm} is also true. Here is a fundamental example of morphism that respects lcm's (\cf{} \cite{C,G} and theorem \ref{morph ad theo} below)~:

\begin{lem}\label{exemple de lcm-morphisme}Let $(J_\a)_{\a \in \tilde I}$ be a family of non-empty spherical subsets of $I$ and assume that, for all $\a,\, \b \in \tilde I$, $\tilde m_{\a,\b} \neq \infty$ implies that $\G_{J_\a\cup J_\b}$ is spherical and $\gras r_{J_\a \cup J_\b} = \prod_{\tilde m_{\a,\b}}(\gras r_{J_\a},\gras r_{J_\b})$. Then the map $\gras s_{\a} \mapsto \gras r_{J_\a}$ extends to a morphism from $B^+_{\tilde \G}$ to $B^+_\G$. Moreover, if for all $\a,\, \b \in \tilde I$, $\tilde m_{\a,\b} = \infty$ implies that $\G_{J_\a\cup J_\b}$ is non-spherical, then this morphism respects lcm's.
\end{lem}

\proof  The first point is clear since the hypothesis implies $\prod_{\tilde m_{\a,\b}}(\gras r_{J_\a},\gras r_{J_\b}) = \prod_{\tilde m_{\a,\b}}(\gras r_{J_\b},\gras r_{J_\a}) $ if $\tilde m_{\a,\b} \neq \infty$. Let us show the second point. We get $\p(\gras s_\a) = \gras r_{J_\a} \neq 1$ since $J_\a$ is non-empty. Moreover, we have the following sequence of equivalences (where the symbol $\vee$ stands for $\vee_L$ or $\vee_R$)~: $\gras s_\a \vee \gras s_\b$ exists in $B^+_{\tilde \G} \Leftrightarrow \tilde m_{\a,\b} \neq \infty \Leftrightarrow \G_{J_\a \cup J_\b}$ is spherical $\Leftrightarrow \gras r_{J_\a \cup J_\b} =  \gras r_{J_\a} \vee  \gras r_{J_\b}$ exists in $B^+_\G$, in which case we get $\p(\gras s_\a \vee \gras s_\b) = \p(\prod_{ \tilde m_{\a,\b}}(\gras s_{\a},\gras s_{\b})) = \prod_{ \tilde m_{\a,\b}}(\gras r_{J_\a},\gras r_{J_\b}) = \gras r_{J_\a \cup J_\b} = \gras r_{J_\a} \vee\gras r_{J_\b} = \p(\gras s_\a) \vee \p(\gras s_\b)$.
\qed\medskip

\subsection{Admissible morphisms, submonoids and subgroups.}\label{Morphismes admissibles entre mono�des d'Artin-Tits.}\mbox{}\medskip

Let $\G = (m_{i,j})_{i,j\in I}$ be a Coxeter matrix. 

The admissibility of a spherical partition $\tilde I$ of $\G$ can naturally be expressed in terms of simple elements in $B^+_\G$. Indeed, if we denote by $\gras {\tilde W}$ the image of the subgroup $\tilde W = \langle r_\a \mid \a \in \tilde I \rangle$ of $W_\G$ in $\gras W_\G \subseteq B^+_\G$, then we get that $\tilde I$ is admissible if and only if, for all $(\gras w, \a) \in \gras {\tilde W}\times \tilde I$, either the products $\gras w \cdot \gras s_i$ are simple for all $i \in \a$, or $\gras w \succcurlyeq \gras s_i$ for all $i \in \a$. In the same way, the compatibility of words on $\tilde I$ is easy to characterize~: 

\begin{lem}\label{compatibilit� vs. r�duit}let $\tilde I$ be a spherical partition of $\G$ and fix $\a_1,\ldots,\,\a_n \in \tilde I$. Then 
\[ 
 \textrm{the word } \textstyle{\prod_{k = 1}^n}\a_k \textrm{  is compatible  } \Leftrightarrow  \textrm{the element } \textstyle{\prod_{k = 1}^n}\gras r_{\a_k} \textrm{  is simple}.
\]
In that case, if $w = \prod_{k = 1}^n r_{\a_k}$ in $\tilde W$, then $\gras w = \prod_{k = 1}^n\gras r_{\a_k}$ in $\gras W_\G$.
\end{lem}
\proof  Set $w = \prod_{k = 1}^n r_{\a_k} = \pi\left(\prod_{k = 1}^n\gras r_{\a_k}\right)$. Assume that $ \prod_{k = 1}^n\a_k$ is compatible, \ie{} $\l(w) = \sum_{k = 1}^n\l(r_{\a_k})$, and fix a reduced representation $R_{\a_k}$ of each $r_{\a_k}$ on $I$. Then the representation $\prod_{k = 1}^n R_{\a_k}$ of $w$ on $I$ is reduced and hence, by definition of $\gras w$, we get $\gras w = \prod_{k = 1}^n\gras r_{\a_k}$ in $\gras W_\G$. Conversely, if the product $\prod_{k = 1}^n\gras r_{\a_k}$ is simple, then $\l(w) = \l\left(\prod_{k = 1}^n\gras r_{\a_k}\right) = \sum_{ k = 1}^n\l(\gras r_{\a_k}) = \sum_{k = 1}^n\l(r_{\a_k})$ (the first and third equalities by definition of $\gras W_\G$, and the second by additivity of the standard length on $B^+_\G$), whence the compatibility of $\prod_{k = 1}^n\a_k$.
\qed\medskip

This lemma allows us to reformulate the characterizations of the admissibility of a $2$-partition of $\G$  (\cf{} lemma \ref{caract ad} above) in terms of simple elements of~$B^+_\G$~:

\begin{lem}\label{caract ad 2}Let $\tilde I = \{\a,\b\}$ be a spherical $2$-partition of $\G$. 
\begin{enumerate}
\item The following conditions are equivalent~: 
\begin{enumerate}
\item $\tilde I$ is an admissible partition of $\G$,
\item for every integer $0 \leqslant n < |r_\a r_\b| + 1$, the two elements $\prod_n(\gras r_\a, \gras r_\b)$ and $\prod_n(\gras r_\b, \gras r_\a)$ of $B^+_\G$ are simple.
\end{enumerate}
\item If $\G$ is spherical, then $|r_\a r_\b| \neq \infty$ and conditions (\ref{caract ad1}) and (\ref{caract ad3}) above are equivalent to the following~:
\begin{enumerate}
 \item the elements $\prod_{|r_\a r_\b| }(\gras r_\a, \gras r_\b)$ and $\prod_{|r_\a r_\b| }(\gras r_\b, \gras r_\a)$ of $B^+_\G$ are simple.
\end{enumerate}
Moreover, we get in that case $\prod_{|r_\a r_\b| }(\gras r_\a, \gras r_\b) = \prod_{|r_\a r_\b| }(\gras r_\b, \gras r_\a) = \gras r_I$.
\end{enumerate}
\end{lem}

We are now able to prove the analogue of theorem \ref{add implique inj} for Artin-Tits monoids and for Artin-Tits groups of spherical type~:

\begin{thm}\label{morph ad theo}Let $\tilde I$ be an admissible partition of $\G$, of type $\tilde \G$. Then~: 
\begin{enumerate}
\item the map $\gras S_{\tilde \G} \to B^+_\G$, $\gras s_\a \mapsto \gras r_\a$, extends to a morphism $\p = \p_{\tilde I} : B^+_{\tilde \G} \to B^+_\G$, 
\item this morphism respects lcm's, hence is injective. 
\end{enumerate}

In particular, if we set $\tilde {\gras S} = \{ \gras r_\a \mid \a \in \tilde I\}$ and denote by $\tilde B^+ = \langle \tilde {\gras S}\rangle^+$ the submonoid of $B^+_\G$ generated by the $\gras r_\a$, $\a \in \tilde I$, then the pair $(\tilde B^+, \tilde {\gras S})$ is (isomorphic to) the positive Artin-Tits system of type $\tilde \G$.
\end{thm}
\proof  We can apply lemma \ref{exemple de lcm-morphisme} to the set $\tilde I$, since it consists of non-empty spherical subsets of $I$, and since we have $|r_\a r_\b| \neq \infty$ if and only if $\G_{\a\cup \b}$ is spherical (by proposition \ref{I spherique ssi tildeI sph�rique}), in which case we get $\prod_{|r_\a r_\b|}(\gras r_\a, \gras r_\b) = \gras r_{\a \cup \b}$ by lemma \ref{caract ad 2}. 
\qed\medskip

The morphism $\p : B^+_{\tilde \G} \hookrightarrow B^+_\G$ of theorem \ref{morph ad theo} clearly extends to a group homomorphism $\p_{\mathrm{gr}} : B_{\tilde \G} \to B_\G$ whose image is the subgroup $\tilde B = \langle \gras r_\a, \, \a\in \tilde I\rangle$ of $B_\G$. When $\tilde \G$ is spherical, the injectivity of $\p$ implies the following~:

\begin{thm}\label{morphismes admissibles et groupes d'Artin}Let $\tilde I$ be an admissible partition of $\G$, of spherical type $\tilde \G$. Then the homomorphism $\p_{\mathrm{gr}} : B_{\tilde \G} \to B_\G$ is injective. In other words, the pair $(\tilde B, \tilde {\gras S})$ is (isomorphic to) the Artin-Tits system of type $\tilde \G$.
\end{thm}
\proof  Since $\tilde \G$ is spherical, every $b \in B_{\tilde \G}$ can be written $b = x^{-1}y$ for $x,\, y \in B^+_{\tilde \G}$ (\cf{} subsection \ref{Generalities on Coxeter and Artin-Tits groups}), and the equality $\p_{\mathrm{gr}}(b) = 1$ hence implies $\p(x) = \p(y)$, whence the result thanks to the injectivity of $\p$.
\qed\medskip

Let us name the objects we have just defined~:

\begin{defn}\label{morphs ad}Let $J \subseteq I$ be a subset of $I$ and let $\tilde J$ be an admissible partition of $\G_J$, of type $\tilde \G$. Let $\tilde {\gras S} = \{\gras s_\a \mid \a \in \tilde J \}$. Then we say that~: 
\begin{enumerate}
\item[$\bullet$] the submonoid $\tilde B^+ = \langle \tilde {\gras S}\rangle^+$ of $B^+_\G$ (resp. the sugroup $\tilde B = \langle \tilde {\gras S} \rangle$ of $B_\G$) is {\em induced\/} by $\tilde J$, or, by abuse of language, is an {\em admissible\/} submonoid (resp. subgroup) of $B^+_\G$ (resp. $B_\G$),
\item[$\bullet$] the morphism $\p = \p_{\tilde J} : B^+_{\tilde \G} \hookrightarrow B^+_\G$ (resp. $\p_{\mathrm{gr}} : B_{\tilde \G} \to B_\G$), which sends each $\gras s_\a \in \gras S_{\tilde \G}$ on $\gras r_\a \in \tilde {\gras S}$, is {\em induced\/} by $\tilde J$, or, by abuse of language, is an {\em admissible\/} morphism.
\end{enumerate}
\end{defn}

\begin{rem}In our definitions, we allow partitions of {\em subsets\/} of $I$. This generalization does not change the conclusions of theorems \ref{morph ad theo} and \ref{morphismes admissibles et groupes d'Artin}, and allows the notion of admissible submonoids, subgroups or morphisms, to comprise the notions of standard parabolic submonoids or subgroups, of submonoids of fixed elements under the action of graph automorphisms and of LCM-homomorphisms of { \cite{C, G}\/} (see theorems \ref{theo sym Crisp} and \ref{correspondance} below). 
\end{rem}

\begin{rem}If the partition $\tilde J$ of $\G_J$ is only supposed to be spherical, then the map $\gras S_{\tilde \G} \to B^+_\G$, $\gras s_\a \mapsto \gras r_\a$, does not necessarily extend to a morphism from $B^+_{\tilde \G}$ to $B^+_\G$~: for example, if $\G =$
\begin{picture}(45,10)(-1,-2)
\put(0,0){\circle{4}}
\put(20,0){\circle*{4}}
\put(40,0){\circle*{4}}
\put(0,3){\scriptsize $1$}
\put(20,3){\scriptsize $2$}
\put(40,3){\scriptsize $3$}
\put(2,0){\line(1,0){18}}
\put(20,0){\line(1,0){20}}
\end{picture}
with $\a = \{1\}$ and $\b = \{2,3\}$, then $|r_\a r_\b| = 3$ but $\gras r_\a \gras r_\b \gras r_\a \neq \gras r_\b\gras r_\a\gras r_\b$ in $B^+_\G$ (look at the standard length).
\end{rem}

\subsection{Admissibility and Artin-Tits monoids or groups in the literature.}\label{Les morphismes admissibles dans la litt�rature.}\mbox{}\medskip

In this subsection, we show how our notions of admissible submonoids, subgroups or morphisms generalize and unify three situations that have been studied earlier.

\subsubsection{Submonoids of fixed points under the action of graph automorphisms.}\label{points fixes sous sigma}\mbox{}\medskip

Here is the analogue of \cite[Cor.~3.5]{H} and \cite[Thm.~1.3]{Mu} (\cf{} subsection \ref{exemples de partitions admissibles} above) for Artin-Tits monoids and for Artin-Tits groups of spherical type. Hence we recover the results \cite[Thm.~9.3]{DP}, \cite[Cor.~4.4]{Mi} and \cite[Lem.~10 and Thm.~11]{C2}.

\begin{thm}\label{theo sym Crisp}Let $\G = (m_{i,j})_{i,j \in I}$ be a Coxeter matrix and $G$ be a subgroup of $\Aut(\G)$. Let $\tilde J$ be the set of all spherical orbits of $I$ under $G$ and let $J \subseteq I$ be their union. Let $\tilde \G$ be the type of the admissible partition $\tilde J$ of $\G_J$, and set $\tilde {\gras S} = \{\gras r_\a \mid \a \in \tilde J\}$, $\tilde B^+ = \langle \tilde {\gras S} \rangle^+$ and $\tilde B = \langle \tilde {\gras S} \rangle$. Then~:
\begin{enumerate}
\item $(B^+_\G)^G = \tilde B^+$ and hence the pair $((B^+_\G)^G,\tilde {\gras S})$ is (isomorphic to) the positive Artin-Tits system of type $\tilde \G$,
\item if $\G$ is spherical, then $(B_\G)^G = \tilde B$ and hence the pair $((B_\G)^G,\tilde {\gras S})$ is (isomorphic to) the Artin-Tits system of type $\tilde \G$.
\end{enumerate}
\end{thm}
\proof  We already know that $\tilde J$ is an admissible partition of $\G_J$ (\cf{} subsection \ref{exemples de partitions admissibles}). Thanks to theorems \ref{morph ad theo} and \ref{morphismes admissibles et groupes d'Artin} above, the only things to prove are $(B^+_\G)^G = \tilde B^+$ and, when $\G$ is spherical, $(B_\G)^G = \tilde B$.

For $\a \in \tilde J$, the group $G$ stabilizes $B^+_\a$ and the induced action respects the standard length, so $G$ fixes $\gras r_\a$ (which is the unique element of maximal standard length in $\gras W_\a$). Hence we get $\tilde B^+ \subseteq (B^+_\G)^G$ and $\tilde B \subseteq (B_\G)^G$.

Let $x$ be an element of $(B^+_\G)^G$ and let us show by induction on $\l(x)$ that $x \in \tilde B^+$. There is nothing to prove if $x = 1$, so assume $x \neq 1$ and consider an element $i \in I$ such that $\gras s_i \preccurlyeq x$. Then, for all $g \in G$, $\gras s_{g(i)} \preccurlyeq x$. This implies that the orbit $\a$ of $i$ under $G$ is spherical and that $\gras r_\a \preccurlyeq x$. So there exists $x' \in B^+_\G$ such that $x = \gras r_\a x'$, and $\l(x')<\l(x)$. By cancellativity in $B^+_\G$, we get $x' \in (B^+_\G)^G$, hence $x' \in \tilde B^+$ by induction, and finally $x \in \tilde B^+$.

Now assume that $\G$ is spherical and fix $b \in (B_\G)^G$. Let $(x,y) \in (B^+_\G)^2$ be the irreducible left form of $b$ (i.e. the unique pair such that $b = x^{-1}y$ and $x\wedge_L y =1$, \cf{} definition \ref{def irred left form} above). Since the action of $G$ on $B^+_\G$ respect divisibility (hence gcd's), we get by unicity that $x,\, y \in (B^+_\G)^G$. The first point then gives $x,\,y \in \tilde B^+$, whence $b \in \tilde B$. 
\qed\medskip

\begin{rem}\label{sym type FC}On the work of Crisp \cite{C2}.
\begin{enumerate}
\item Our proof of theorem \ref{theo sym Crisp} is very similar to those of \cite[Lem.~10 and Thm.~11]{C2}, and indeed, the results \cite[Lem.~6]{C2}, \cite{C2'} and lemma \ref{LCM-donn�e et ordre} below show that the Coxeter matrix $(m_{BC})_{B,C \in {\bf S}}$ constructed by Crisp in  \cite{C2, C2'} is precisely our matrix $\tilde \G$. 
\item Crisp actually established the second point of theorem \ref{theo sym Crisp} for a wider class of Coxeter graphs than the spherical ones, namely the type FC ones, \ie{} the finite Coxeter graphs for which every complete subgraph with no infinite label is spherical \cite[Thm.~4]{C2}.
\end{enumerate}
\end{rem}

\subsubsection{LCM-homomorphisms.}\label{Morphismes admissibles et LCM-homomorphismes.}\mbox{}\medskip

We recall in definition \ref{LCM-donn�es} below the notion of {\em LCM-homomorphisms\/} as defined in \cite[Def.~2.1]{G}, which generalizes the one of \cite[Def.~2.1]{C} by allowing finite Coxeter graphs with infinite labels. We adapt these definitions to our settings by defining the notion of  {\em LCM-partitions\/} of a Coxeter graph, which will turn out to be nothing else but special cases of admissible partitions (\cf{} proposition~\ref{correspondance} below). We do not suppose that the Coxeter graphs involved are finite.

\begin{defn}\label{LCM-partition}Let $\G = (m_{i,j})_{i,j\in I}$ be a Coxeter matrix and let $\tilde I$ be a spherical partition of $\G$. Let $\Omega = (n_{\a,\b})_{\a,\b \in \tilde I}$ be a Coxeter matrix over $\tilde I$. We say that $\tilde I$ is an {\em LCM-partition\/} of $\G$, of {\em type\/} $\Omega$, if, for each pair $(\a, \b) \in \tilde I ^2$, we have the following alternative~: 
\begin{enumerate}
\item[(Fi)] $n_{\a,\b} \neq \infty$, $\G_{\a \cup \b}$ is spherical and $\gras r_{\a \cup \b} = \prod_{n_{\a,\b}}(\gras r_\a,\gras r_\b)$,
\item[(In)] $n_{\a,\b} = \infty$ and for all $i \in \a$, $\G_{\{i\}\cup \b}$ is non-spherical. 
\end{enumerate}
\end{defn}

\begin{defn}[{\cite[Defs.~2.1]{C, G}}]\label{LCM-donn�es}Let $\G = (m_{i,j})_{i,j\in I}$ be a Coxeter matrix. Let $J \subseteq I$ be a subset of $I$ and let $\tilde J$ be an LCM-partition of $\G_J$, of type $\Omega = (n_{\a,\b})_{\a,\b \in \tilde J}$. Lemma \ref{exemple de lcm-morphisme} above shows that the map ${\gras S}_\Omega \to B^+_\G$, $\gras s_\a \mapsto \gras r_\a$, extends to a morphism that respects lcm's from $B^+_\Omega$ to $B^+_\G$, which we call, after \cite{C,G}, an {\em LCM-homomorphism}.
\end{defn}

Let $\G = (m_{i,j})_{i,j\in I}$ be a Coxeter matrix, and let $\tilde I$ be an LCM-partition of $\G$, of type $\Omega = (n_{\a,\b})_{\a,\b \in \tilde I}$. We show in proposition \ref{correspondance} below that $\tilde I$ is an admissible partition of $\G$, and that its type (as an LCM-partition) $\Omega$ is necessarily its type (as a spherical partition) $\tilde \G = (|r_\a r_\b|)_{\a,\b \in \tilde I}$.

\begin{lem}\label{LCM-donn�e et ordre}Let $\G = (m_{i,j})_{i,j \in I}$ be a Coxeter matrix and let $\a$ and $\b$ be two spherical subsets of $I$.
\begin{enumerate}
\item \label{LCM-donn�e et ordre1} If $\G_{\a \cup \b}$ is spherical and if there exists an integer $n \in \Nset$ such that $\gras r_{\a \cup \b} = \prod_{n}(\gras r_\a,\gras r_\b) = \prod_{n}(\gras r_\b,\gras r_\a)$, then $n = |r_\a r_\b|$.
\item \label{LCM-donn�e et ordre2} If, for all $n \in \Nset$, the product $\prod_{n}(\gras r_\a,\gras r_\b)$ is simple, then $|r_\a r_\b | = \infty$ and $\G_{\a \cup \b}$ is non-spherical.
\end{enumerate} 
\end{lem}
\proof  Under the hypothesis of assertion {(\ref{LCM-donn�e et ordre1})}, we get $(r_\a r_\b)^{n} = \prod_{2 n}(r_\a,r_\b) = (r_{\a\cup \b})^2 = 1$ in $W_\G$, hence $| r_\a r_\b|$ divides $n$. If $| r_\a r_\b | < n$, then we can replace a factor $\prod_{|r_\a r_\b|}(r_\a,r_\b)$ of $\prod_{n}(r_\a,r_\b)$ by $\prod_{|r_\a r_\b|}(r_\b,r_\a)$ and then simplify $2|r_\a r_\b|$ terms, whence $\l\left(\prod_{n}(r_\a,r_\b)\right) < \sum_n(\l(r_\a) ,\l(r_\b)) = \sum_n(\l(\gras r_\a) ,\l(\gras r_\b)) = \l(\prod_{n}(\gras r_\a,\gras r_\b))$, and a contradiction since $\prod_{n}(\gras r_\a,\gras r_\b)$ is simple. Under the hypothesis of assertion {(\ref{LCM-donn�e et ordre2})}, the dihedral group $\langle r_\a,r_\b\rangle$, which is included in $W_{\a\cup\b}$, is infinite, hence $|r_\a r_\b | = \infty$ and $\G_{\a \cup \b}$ is non-spherical.
\qed\medskip

\begin{thm}\label{correspondance}Let $\G = (m_{i,j})_{i,j\in I}$ be a Coxeter matrix, and let $\tilde I$ be an LCM-partition of $\G$, of type $\Omega = (n_{\a,\b})_{\a,\b \in \tilde I}$. Then $\tilde I$ is an admissible partition of $\G$, and $\Omega = \tilde \G = (|r_\a r_\b|)_{\a,\b \in \tilde I}$.
\end{thm}
\proof  A consequence of \cite[Lem.~2.5]{G} is that, if $n_{\a,\b} = \infty$, then for all $n \in \Nset$, the product $\prod_{n}(\gras r_{\a},\gras r_{\b})$ is simple. Lemma \ref{LCM-donn�e et ordre} then shows that $\Omega = \tilde \G$ and the characterizations of lemma \ref{caract ad 2} show that for all $\a,\,\b \in \tilde I$, $\{\a,\b\}$ is an admissible partition of $\G_{\a \cup \b}$. We conclude that $\tilde I$ is an admissible partition of $\G$ thanks to theorem \ref{admissibilit� locale}. 
\qed\medskip 

So, as announced, an LCM-partition is an admissible partition (and hence an LCM-homomorphism is an admissible morphism)~; the converse is false in general (\cf{} example \ref{sym�trie de A3tilde}, remark \ref{eclatement pas LCM} and example \ref{L3b mais pas L3c dans un type FC} below), but is true for example if~:
\begin{enumerate}
\item the matrix $\tilde \G$ has no infinite coefficient,
\item the matrix $\G$ is {\em right angled}, \ie{} $m_{i,j} \in \{1,2,\infty\}$ for all $i,\,j \in I$ (to see this, use \cite[Lem.~2.5.15]{Mu2}, recalled in proposition~\ref{2.5.15} below),
\item the matrix $\G$ is of type FC (this notion is defined in remark \ref{sym type FC}) and $\tilde I$ is the set of orbits of $I$ under the action of a subgroup of $\Aut(\G)$.
\end{enumerate}

\begin{exmp}\label{sym�trie de A3tilde}Consider the Coxeter graph $\G$ of affine type $\tilde A_3$, and its $2$-partition formed by pairs of opposite vertices~: 
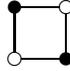
\begin{center}
\begin{picture}(20,30)(0,-5)\setlength{\unitlength}{1.3pt}
\put(0,0){\circle{4}}
\put(15,0){\circle*{4}}

\put(0,15){\circle*{4}}
\put(15,15){\circle{4}}

\put(2,0){\line(1,0){12}}
\put(0,2){\line(0,1){12}}
\put(15,0){\line(0,1){13}}
\put(0,15){\line(1,0){13}}
\end{picture}
\end{center}
This spherical $2$-partition is admissible since it is the set of orbits of $\G$ under the action of the "central symmetry", and its type is $\tilde \G = I_2(\infty)$ since $\G$ is non-spherical. It is not an LCM-partition (condition {(In)\/} of definition~{\ref{LCM-partition}\/} is not satisfied)~: indeed, if $i$ is one of the vertex of $\G$ and if $\b$ is the orbit that does not contain $i$, then $\G_{\{i\}\cup\b}$ is of spherical type $A_3$. 
\end{exmp}

\begin{rem}
 The results { \cite[Prop. 2.3]{C}\/} and { \cite[Cor. 2.7]{G}\/} on the injectivity of the morphism between Coxeter groups induced by an LCM-homomorphism now appear as special cases of { \cite[Thm. 1.1]{Mu}\/} (recalled in theorem { \ref{add implique inj}} above). In fact, one can check that the proof of { \cite[Cor. 2.7]{G}\/} works for general admissible partitions and hence gives a new proof of { \cite[Thm. 1.1]{Mu}}.
\end{rem}

\subsubsection{The bursts of a Coxeter graph.} \label{The bursts of a Coxeter graph.}\mbox{}\medskip

We recall here a construction of M\"uhlherr \cite[section~2.6]{Mu2}, a quasi-identical version of which has independently been obtained by Crisp and Paris for Coxeter graphs with no infinite label \cite[section~6]{CP}, and by Paris in general \cite[section~5]{P}. The differences between the two approaches rely essentially in the choice of the integer $N$ in definition \ref{def eclatement} below.

Let $\delta : \Nset_{\geqslant 2} \cup\{\infty\} \to \Nset_{\geqslant 1}, m \mapsto \begin{cases} m-1 &\textrm{  if $m$ is even},\\ \dfrac{m-1}{2} & \textrm{  if $m$ is odd},\\  2 & \textrm{  if $m = \infty$}. \end{cases}$

\begin{defn}[{\cite[section~2.6]{Mu2}}]\label{def eclatement}Suppose that $\G = (m_{i,j})_{i,j \in I}$ is a Coxeter matrix such that the subset $\{m_{i,j} \mid i,\, j\in I\}$ of $\Nset \cup \{\infty\}$ is finite. Set $N_0 = \lcm \{\delta(m_{i,j}) \mid i,\, j \in I, \, i\neq j\}$ and let $N$ be a multiple of $N_0$. A {\em $N$-burst}, or simply a {\em burst}, of $\G$ is a Coxeter graph $\widehat \G$ with vertex set the disjoint union $\widehat I = \bigsqcup_{i \in I} T(i)$ of sets $T(i) = \{i^{(1)},\ldots , i^{(N)}\}$  of cardinality $N$, and with edges displayed as follows~: 
\begin{enumerate}
\item there is no edge between two elements of a same $T(i)$,
\item if $m_{i,j} \in \Nset_{\geqslant 2}$ is even, the graph $\widehat \G_{T(i) \sqcup T(j)}$ is the disjoint union of  $\frac{N}{\delta(m_{i,j})}$ copies of the following graph~:
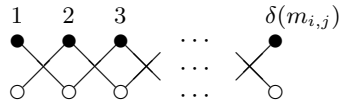
\begin{center}
\begin{picture}(125,48)(0,-10)\setlength{\unitlength}{1.3pt}
\put(0,0){\circle{4}}
\put(15,0){\circle{4}}
\put(30,0){\circle{4}}
\put(75,0){\circle{4}}

\put(0,15){\circle*{4}}\put(-2,20){\small $1$}
\put(15,15){\circle*{4}}\put(13,20){\small $2$}
\put(30,15){\circle*{4}}\put(28,20){\small $3$}
\put(75,15){\circle*{4}}\put(72,20){\small $\delta(m_{i,j})$}

\put(1.5,1.5){\line(1,1){14}}
\put(16.5,1.5){\line(1,1){14}}
\put(31.5,1.5){\line(1,1){10}}
\put(13.5,1.5){\line(-1,1){14}}
\put(28.5,1.5){\line(-1,1){14}}
\put(31.5,13.5){\line(1,-1){10}}
\put(73.5,13.5){\line(-1,-1){10}}
\put(73.5,1.5){\line(-1,1){10}}

\put(47,-3){$\cdots$}
\put(47,5){$\cdots$}
\put(47,13){$\cdots$}
\end{picture}
\end{center}
where the vertices $\bullet$ constitute $T(i)$ and the vertices $\circ$ constitute $T(j)$,
\item if $m_{i,j} \in \Nset_{\geqslant 3}$ is odd, the graph $\widehat \G_{T(i) \sqcup T(j)}$ is the disjoint union of  $\frac{N}{\delta(m_{i,j})}$ copies of the following graph~:
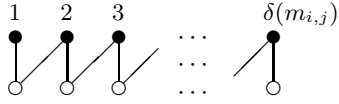
\begin{center}
\begin{picture}(125,48)(0,-10)\setlength{\unitlength}{1.3pt}
\put(0,0){\circle{4}}
\put(15,0){\circle{4}}
\put(30,0){\circle{4}}
\put(75,0){\circle{4}}

\put(0,15){\circle*{4}}\put(-2,20){\small $1$}
\put(15,15){\circle*{4}}\put(13,20){\small $2$}
\put(30,15){\circle*{4}}\put(28,20){\small $3$}
\put(75,15){\circle*{4}}\put(72,20){\small $\delta(m_{i,j})$}

\put(1.5,1.5){\line(1,1){14}}
\put(16.5,1.5){\line(1,1){14}}
\put(31.5,1.5){\line(1,1){10}}
\put(0,2){\line(0,1){13}}
\put(15,2){\line(0,1){13}}
\put(30,2){\line(0,1){13}}
\put(73.5,13.5){\line(-1,-1){10}}
\put(75,2){\line(0,1){13}}

\put(47,-3){$\cdots$}
\put(47,5){$\cdots$}
\put(47,13){$\cdots$}
\end{picture}
\end{center}
where the vertices $\bullet$ constitute $T(i)$ and the vertices $\circ$ constitute $T(j)$,
\item if $m_{i,j}  = \infty$, the graph $\widehat \G_{T(i) \sqcup T(j)}$ is the disjoint union of  $\frac{N}{\delta(m_{i,j})}$ copies of the following graph~:
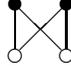
\begin{center}
\begin{picture}(125,30)(0,-6)\setlength{\unitlength}{1.3pt}
\put(30,0){\circle{4}}
\put(45,0){\circle{4}}
\put(30,15){\circle*{4}}
\put(45,15){\circle*{4}}

\put(31.5,1.5){\line(1,1){14}}
\put(43.5,1.5){\line(-1,1){14}}
\put(45,2){\line(0,1){13}}
\put(30,2){\line(0,1){13}}
\end{picture}
\end{center}
where the vertices $\bullet$ constitute $T(i)$ and the vertices $\circ$ constitute $T(j)$.
\end{enumerate}
\end{defn}

\begin{thm}[{\cite[Thm.~2.6.1 and its proof]{Mu2}}]\label{�clatement}Let $\G = (m_{i,j})_{i,j \in I}$ be a Coxeter matrix with $\{m_{i,j} \mid i,\, j\in I\}$ finite, and let $\widehat \G$ be a $N$-burst of $\G$. Then the partition $\{T(i) \mid i \in I\}$ of $\widehat I$ is an admissible partition of $\widehat\G$, of type (isomorphic to) $\G$.
\end{thm}
\proof  It is enough to check that, for all $i,\, j \in I$, $i\neq j$, $\{T(i), T(j)\}$ is an admissible partition of $\widehat \G_{T(i) \sqcup T(j)}$, of type $I_2(m_{i,j})$ (with $I_2(2) = A_1\times A_1$).

If $m_{i,j} = 2$, then there is no edge between a vertex of $T(i)$ and a vertex of $T(j)$. If $m_{i,j} \in \Nset_{\geqslant 3}$, then the graph $\widehat \G_{T(i) \sqcup T(j)}$ is the disjoint union of $\frac{2N}{m_{i,j}-1}$ copies of the spherical Coxeter graph of type $A_{m_{i,j}-1}$, and the partition $\{T(i), T(j)\}$ induces on each of these connected components the {\em bipartite\/} partition of $A_{m_{i,j}-1}$. If $m_{i,j} = \infty$, then the graph $\widehat \G_{T(i) \sqcup T(j)}$ is the disjoint union of $\frac{N}{2}$ copies of the affine Coxeter graph of type $\tilde A_3$, and the partition $\{T(i), T(j)\}$ induces on each of these connected components the partition of $\tilde A_3$ described in example \ref{sym�trie de A3tilde} above. We conclude by applying results of \cite[section~2.5]{Mu2} recalled in propositions \ref{admissibili� et tilde m=2}, \ref{cas des �clatements n} and \ref{ad des partitions bipartites} below (note that we really need our stronger version, prop.~\ref{cas des �clatements n}, of \cite[Lem.~2.5.4]{Mu2} when $m_{i,j} = \infty$). \qed\medskip

\begin{exmp}\label{eclatement H3 H4}If $\G$ is of type $H_3$ (resp. $H_4$), then $N_0 = 2$ and every $2$-burst $\widehat \G$ of $\G$ is of type $D_6$ (resp. $E_8$). We thus recover the figures of example { \ref{H3 dans D6}}.
\end{exmp}

\begin{rem}\label{eclatement pas LCM}When $\G$ has an infinite coefficient, then $\{T(i) \mid i \in I\}$ is not an LCM-partition of $\widehat \G$ (condition { (In)\/} of definition~{ \ref{LCM-partition}\/} is not satisfied)~: indeed, if $m_{i,j} = \infty$, then for $i^{(k)} \in T(i)$, we get that the graph $\widehat\G_{\{i^{(k)}\}\cup T(j)}$ is the disjoint union of $N-2$ connected components of type $A_1$ and one connected component of type $A_3$, hence is spherical. 
\end{rem}


\subsection{Some properties of admissible morphisms.}\label{Quelques propri�t�s des morphismes admissibles.}\mbox{}\medskip

In this subsection, we show that some properties established in \cite{C, C2, G} for their special cases of admissible morphisms are in fact satisfied by all admissible morphisms.

\subsubsection{Respect of the combinatorics.}\mbox{}\medskip

Let $\G = (m_{i,j})_{i,j \in I}$ be a Coxeter matrix, $J \subseteq I$ be a subset of $I$, and $\tilde J$ be an admissible partition of $\G_J$ of type $\tilde \G$. We consider the admissible morphism $\p : B^+_{\tilde \G} \hookrightarrow B^+_\G$ induced by $\tilde J$, and we denote by $\gras {\tilde W}$ the image of the subgroup $\tilde W = \langle r_\a ,\, \a \in \tilde J \rangle$ of $W_\G$ in  $\gras W_\G \subseteq B^+_\G$.

We know that $\p$ respects lcm's and divisibility, in the sense of theorem \ref{props morphisme ppcm} above. The following lemma establishes that $\p$ respects the notions of simple elements in $B^+_{\tilde \G}$ and in $\tilde B^+_\G$~; it is a generalization of the well-known analogous result for the standard parabolic subgroups, and of \cite[Lem.~15]{C2}, \cite[Lem.~2.2]{C} and \cite[Prop.~2.6]{G}.

\begin{lem}\label{[C, 2.2], [C2, 15] et [G, 2.6]}With the above notations, we get $\p(\gras W_{\tilde \G}) = \tilde B^+ \cap \gras W_\G = \gras {\tilde W}$. Moreover, if $\tilde \G$ (or equivalently $\G_J$) is spherical, then $\p(\gras r_{\tilde J}) = \gras r_J$.
\end{lem}
\proof  This is a direct consequence of proposition~\ref{I spherique ssi tildeI sph�rique} and lemma~\ref{compatibilit� vs. r�duit}.\qed\medskip

Let us mention two consequences of that result, given by \cite[Thm.~2.10 and Cor.~2.11]{G}, which apply to our settings~; note however that for the proofs of \cite[Lem.~2.9 and Thm.~2.10]{G} to be correct, we have to add to their hypothesis the following condition, which is satisfied by any admissible morphism~: $\mathrm{Im}(\p) \subseteq B^+_{\bigcup_{\a \in \tilde J}p(\a)}$, where $p(\a) = \{i \in I \mid \gras s_i \preccurlyeq \p(\gras s_\a) \} = \{i \in I \mid \p(\gras s_\a) \succcurlyeq \gras s_i \}$. 

\begin{prop}[{\cite[Thm.~2.10]{G}}]\label{respect des formes normales}Let $\p$ be as above. Then~: 
\begin{enumerate}
\item the morphism $\p$ respects (left and right) normal forms, \ie{} if $(x_1,\ldots,\, x_n)$ is the left (resp. right) normal form of a non trivial element $x \in B^+_{\tilde \G}$, then $(\p(x_1),\ldots,\, \p(x_n))$ is the left (resp. right) normal form of $\p(x) \in B^+_\G$,
\item the morphism $\p$ respects (left and right) gcd's, \ie{} for all $(x, y) \in (B^+_{\tilde \G})^2$, we get $\p(x \wedge_L y) = \p(x) \wedge_L \p(y)$ and $\p(x \wedge_R y) = \p(x) \wedge_R \p(y)$.
\end{enumerate}
\end{prop}

\begin{cor}[{\cite[Cor.~2.11]{G}}]Assume that $\G$ and $\tilde \G$ are spherical. Then the morphism $\p_{\mathrm{gr}} : B_{\tilde \G} \hookrightarrow B_\G$ respects (left and right) irreducible fractions, \ie{} if $(x,y) \in (B^+_{\tilde \G})^2$ is the left (resp. right) irreducible form of an element $g \in B_{\tilde \G}$, then $(\p(x),\p(y))$ is the left (resp. right) irreducible form of $\p_{\mathrm{gr}}(g) \in B_{\G}$.
\end{cor}

\subsubsection{Composition of admissible morphisms.}\mbox{}\medskip

In proposition \ref{[M�2, 2.5.6]} below, we recall the result \cite[Lem.~2.5.6]{Mu2} on {\em admissible partitions of an admissible partition}. This result implies that the class of admissible morphisms is closed by composition (see corollary \ref{composotion de morphismes ad} below) and offers a criterion to test the admissibility of some spherical partitions, which we use in example \ref{L3b mais pas L3c dans un type FC} below and further in section \ref{Classification}.

\begin{prop}[{\cite[Lem.~2.5.6]{Mu2}}]\label{[M�2, 2.5.6]}Let $\G = (m_{i,j})_{i,j \in I}$ be a Coxeter matrix and let $I'$ be an admissible partition of $\G$, of type $\G'$. Let $I''$ be a spherical partition of $\G'$, of type $\G''$. Set $\overline \Phi = \bigcup_{\a \in \Phi}\a$ for $\Phi \in I''$ and $\overline I = \{\overline \Phi \mid \Phi\in I''\}$. Then $\overline I$ is a spherical partition of $\G$, of type (isomorphic to) $\G''$, and $\overline I$ is admissible if and only if $I''$ is admissible.
\end{prop}

The following result has been established for the LCM-homomorphisms of \cite{C} (\cf{} \cite[page~134]{C}). It can be shown that it is not true for the LCM-homomorphisms of \cite{G}.

\begin{cor}\label{composotion de morphismes ad}The composition of two admissible morphisms is an admissible morphism. 
\end{cor}
\proof Let $\G$, $\G'$ and $\G''$ be three Coxeter matrices and let $\p : B^+_{\G'} \to B^+_\G$ and $\p' : B^+_{\G''} \to B^+_{\G'}$ be two admissible morphisms. In other words, $\G'$ is the type of an  admissible partition $J'$ of $J \subseteq I$, and $\G''$ is the type of an admissible partition $K''$ of $K' \subseteq J'$. But $K'$ is then an admissible partition of $K = \bigcup_{\a \in K'}\a \subseteq J$ (\cf{} theorem \ref{admissibilit� locale}), and proposition \ref{[M�2, 2.5.6]} tells us that $\overline K = \{\overline \Phi \mid \Phi \in K''\}$  is an admissible partition of $K$. Moreover we get $\p \circ \p'(\gras s_\Phi) = \p(\gras r_\Phi) = \p(\lcm\{\gras s_\a \mid \a \in \Phi\}) = \lcm\{\p(\gras s_\a) \mid \a \in \Phi\} = \lcm \{\gras r_\a \mid \a \in \Phi\} = \gras r_{\overline \Phi}$ for every $\Phi \in K''$ (we use proposition~\ref{props morphisme ppcm} for the third equality). Hence $\p \circ \p'$ is the admissible morphism induced by the admissible partition $\overline K$ of $K$.
\qed\medskip

\begin{exmp}\label{L3b mais pas L3c dans un type FC}Consider the two following Coxeter graphs, where $m \in \Nset_{\geqslant 3}$~: 
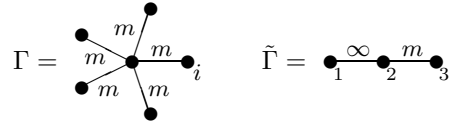
\begin{center}
\begin{picture}(120,40)(-5,-5)
\put(-30,12){$\G = $}
\put(-4,5){\circle*{5.2}}
\put(22,-5){\circle*{5.2}}
\put(15,15){\circle*{5.2}}
\put(-4,25){\circle*{5.2}}
\put(22,35){\circle*{5.2}}
\put(36,15){\circle*{5.2}}\put(38,8){\small $i$}
\put(2,2){\small $m$}
\put(21,2){\small $m$}
\put(22,17){\small $m$}
\put(-3,14){\small $m$}
\put(8,25){\small $m$}
\put(15,15){\line(1,3){7}}
\put(15,15){\line(1,0){21}}
\put(15,15){\line(-2,1){18}}
\put(15,15){\line(-2,-1){18}}
\put(15,15){\line(1,-3){7}}

\put(64,12){$\tilde \G =$}
\put(90,15){\circle*{5.2}}
\put(110,15){\circle*{5.2}}
\put(130,15){\circle*{5.2}}
\put(96,17){\small $\infty$}
\put(117,17){\small $m$}
\put(90,15){\line(1,0){20}}
\put(110,15){\line(1,0){20}}
\put(91,8){\scriptsize $1$}
\put(111,8){\scriptsize $2$}
\put(131,8){\scriptsize $3$}
\end{picture}
\end{center}
The graph $\tilde \G$ (which is of type FC) is the type of the admissible partition of $\G$ composed of orbits of $\G$ under the action of the automorphisms of $\G$ that fix the vertex $i$. Proposition { \ref{[M�2, 2.5.6]}\/} then implies that the spherical partition $\{\{1,3\},\{2\}\}$ of $\tilde \G$ is admissible since it "lifts" to the admissible partition of $\G$ composed of orbits of $\G$ under the action of the whole group $\Aut(\G)$. This admissible $2$-partition of $\tilde \G$ is of type $I_2(\infty)$ (since $\tilde \G$ is not spherical) and is not an LCM-partition (condition { (In)\/} of definition { \ref{LCM-partition}\/} is not satisfied) since $\tilde \G_{\{2,3\}}$ is spherical. 
\end{exmp}

\subsubsection{Geometrical point of view.}\mbox{}\medskip

In \cite[section~3]{C} (resp. in \cite[section~5]{C2} and in \cite[section~3.2]{G}), the authors gave a geometrical interpretation of their special case of admissible morphism between Artin-Tits groups in terms of a map between the associated Salvetti complexes (resp. modified Deligne complexes). One can check that these constructions are still valid for general admissible morphisms.

However, Godelle's proof of the injectivity of LCM-homomorphisms between type FC Artin-Tits groups --- more precisely the proof of \cite[Prop.~3.7]{G} --- does not work for an admissible morphism between type FC Artin-Tits groups that is not an LCM-homomorphism (and such a morphism exists, \cf{} example~\ref{L3b mais pas L3c dans un type FC}). I do not know whether such a morphism is injective or not.

\section{Classification.}\label{Classification}\mbox{}\medskip

The aim of this section is to complete the classification of admissible partitions whose type has no infinite label, began in \cite[section~2.5]{Mu2}. Thanks to our results of subsection \ref{Morphismes admissibles et LCM-homomorphismes.} above, this will in particular give us the classification of LCM-homomorphisms of \cite{C}. 

The results \cite[Thm.~1.2]{Mu} and \cite[Lem.~2.5.5]{Mu2} (\cf{} theorem \ref{admissibilit� locale} and proposition \ref{I spherique ssi tildeI sph�rique} above) reduce this classification  to the classification of admissible $2$-partitions of spherical Coxeter graphs. In subsection \ref{Partitions admissibles et r�ductibilit�.}, we deal with the case $|r_\a r_\b| = 2$ and then recall some results of \cite[section~2.5]{Mu2} which allow to reduce again the problem into the classification of admissible $2$-partitions of irreducible spherical Coxeter graphs.

In subsection \ref{parit� dans An}, we recall the classification of admissible $2$-partitions of Coxeter graphs of types $A_n$, $B_n$ and $D_n$, obtained by M\"uhlherr in \cite[section~2.5]{Mu2}, and complete it by examining the exceptional cases. 

Finally, in subsection \ref{Pliages.}, we compare this classification with the notion of {\em foldings\/} of a Coxeter graph, defined by Crisp in \cite[Def.~4.1]{C} in order to provide examples of LCM-homomorphisms and to begin their classification. This leads us to a generalization (and simplification) of the notion of foldings, which becomes equivalent to the notion of admissible partitions, and allows us to complete the list of cases of the original definition \cite[Def.~4.1]{C}.

\subsection{Admissibility and reducibility.}\label{Partitions admissibles et r�ductibilit�.}\mbox{}\medskip

Let $\G = (m_{i,j})_{i,j\in I}$ be a Coxeter matrix.

Using Tits' solution of the word problem \cite[Thm.~3]{T}, one obtains the following result, where the {\em support\/} of $w \in W_\G$ --- denoted by $\Supp(w)$ --- is the set of letters of any reduced representation of $w$ on $I$ (this set does not depend of the choice of the reduced representation of $w$ since two such words only differ from a finite sequence of braid relations $\prod_{m_{i,j}}(i,j) \leadsto \prod_{m_{i,j}}(j,i)$ with $i,\, j \in I$ such that $i \neq j$ and $m_{i,j}\neq \infty$, which do not change the set of letters involved).

\begin{lem}\label{commutation et supports disjoints}Let $v,\, w \in W_\G$ such that $\Supp(v) \cap \Supp(w) = \emptyset$. Then~:  
\begin{enumerate}
\item \label{commutation et supports disjoints2} $\l(vw) = \l(v) + \l(w)$.
\item \label{commutation et supports disjoints1} $vw = wv \Longleftrightarrow \forall \,(i,j) \in \Supp(v) \times \Supp(w), \ m_{i,j} = 2$,
\end{enumerate}
\end{lem}

We can now deal with the case of the admissible $2$-partitions $\{\a,\b\}$ of $\G$ with $|r_\a r_\b| = 2$~:

\begin{prop}\label{admissibili� et tilde m=2}Let $\tilde I = \{\a,\b\}$ be a spherical $2$-partition of $\G$. Then we have $|r_\a r_\b| = 2 \Leftrightarrow \G = \G_\a \times \G_\b$. In that case, $\tilde I$ is an admissible partition of $\G$.
\end{prop}
\proof If $\G = \G_\a \times \G_\b$, then we obviously have $|r_\a r_\b| = 2$. If $|r_\a r_\b| = 2$, then $r_\a r_\b = r_\b r_\a$ and, by the previous lemma, we get that $\G = \G_\a \times \G_\b$ and $\l(r_\a r_\b) = \l(r_\a) + \l(r_\b)$. The result \cite[Lem.~3.3]{Mu} (\cf{} lemmas \ref{caract ad} or \ref{caract ad 2} above) then implies that $\tilde I$ is an admissible partition of $\G$. 
\qed\medskip

We will need the following proposition to limit the "forms" that an admissible $2$-partition $\{\a,\b\}$ of $\G$ can have when $|r_\a r_\b| \geqslant 3$. For convenience, we sketch the proof of M\"uhlherr below. 

\begin{prop}[{\cite[Lem.~2.5.15]{Mu2}}]\label{2.5.15}Let $\tilde I = \{\a,\b\}$ be an admissible $2$-partition of $\G$. Assume that there exists $i_0 \in \a$ such that $m_{i_0,j} = 2$ for all $j \in \b$. Then $\G = \G_\a \times \G_\b$ (and hence $|r_\a r_\b| = 2$).
\end{prop}
\proof Since $r_\b s_{i_0} = s_{i_0}r_\b$, we get by lemma \ref{commutation et supports disjoints} (first assertion) that $\l(r_\a r_\b s_{i_0}) = \l(r_\a s_{i_0} r_\b) = \l(r_\a s_{i_0}) + \l(r_\b) = \l(r_\a) -1 + \l(r_\b) = \l(r_\a r_\b) -1$, \ie{} $i_0 \in I^-(r_\a r_\b)$. Since $\tilde I$ is admissible, we then have $\a \subseteq I^-(r_\a r_\b)$, whence $I = \a \cup \b \subseteq I^-(r_\a r_\b)$ and $r_\a r_\b = r_I = r_\b r_\a$. We conclude by lemma \ref{commutation et supports disjoints} (second assertion). \qed\medskip

The following proposition allows us to reduce our classification problem to the irreducible case. It is given in \cite[Lem.~2.5.4]{Mu2} for {\em spherical\/} Coxeter graphs $\G_1, \ldots, \G_n$, but in order to complete the proof of theorem \ref{�clatement} above, we need it for general Coxeter graphs. So we prove it below in this more general context, using our characterizations of the admissibility of a $2$-partition of $\G$ in terms of simple elements in $B^+_\G$ (\cf{} lemma \ref{caract ad 2}).

\begin{prop}[{\cite[Lem.~2.5.4]{Mu2}}]\label{cas des �clatements n}Assume that $\G = \G_1\times \cdots \times \G_n$. For $1\leqslant k \leqslant n$, let $\{\a_k,\b_k\}$ be a spherical $2$-partition of $\G_k$ and set $\a = \a_1\sqcup \cdots \sqcup \a_n$ and $\b = \b_1\sqcup \cdots \sqcup \b_n$. Then $\{\a,\b\}$ is a spherical $2$-partition of $\G$ with $r_\a = r_{\a_1}\cdots r_{\a_n}$, $r_\b = r_{\b_1}\cdots r_{\b_n}$ and $| r_\a r_\b | = \lcm\{| r_{\a_k} r_{\b_k} | \mid 1\leqslant k \leqslant n\}$. Moreover, the following conditions are equivalent~:
\begin{enumerate}
\item \label{cas des �clatements n1} $\{\a,\b\}$ is an admissible partition of $\G$,
\item \label{cas des �clatements n2} $\{\a_k,\b_k\}$ is an admissible partition of $\G_k$ for $1\leqslant k \leqslant n$, and $| r_{\a_1} r_{\b_1} | = | r_{\a_2} r_{\b_2} | = \cdots = | r_{\a_n} r_{\b_n} |$.
\end{enumerate}
In that case, we get $| r_\a r_\b | = | r_{\a_1} r_{\b_1} | = | r_{\a_2} r_{\b_2} | = \cdots = | r_{\a_n} r_{\b_n} |$. 
\end{prop}
\proof It is enough to prove the result for $n = 2$. The firsts observations are clear (if needed with the help of lemma \ref{commutation et supports disjoints}). Note that, thanks to lemma \ref{commutation et supports disjoints}, we get that $\gras W_\G \approx \gras W_{\G_1} \times \gras W_{\G_2}$. For example, we have $\gras r_\a = \gras r_{\a_1}\gras r_{\a_2}$ and $\gras r_\b = \gras r_{\b_1}\gras r_{\b_2}$. Hence,  for all $m \in \Nset$, we get $\prod_{m}(\gras r_\a,\gras r_\b) = \prod_{m}(\gras r_{\a_1},\gras r_{\b_1})\prod_{m}(\gras r_{\a_2},\gras r_{\b_2})$ in $B^+_\G \approx B^+_{\G_1} \times B^+_{\G_2}$. 

Suppose { (\ref{cas des �clatements n2})\/} and let us show { (\ref{cas des �clatements n1})}. We get $| r_\a r_\b | = | r_{\a_1} r_{\b_1} | = | r_{\a_2} r_{\b_2} |$. For $k = 1,2$, lemma \ref{caract ad 2} gives us that $\prod_{m}(\gras r_{\a_k},\gras r_{\b_k})$ and $\prod_{m}(\gras r_{\b_k},\gras r_{\a_k})$ are simple for all $0 \leqslant m < | r_{\a_k} r_{\b_k} | +1$. Then $\prod_{m}(\gras r_{\a},\gras r_{\b})$ and $\prod_{m}(\gras r_{\b},\gras r_{\a})$ are simple for all $0 \leqslant m < | r_{\a} r_{\b} | +1$ and we are done by applying lemma \ref{caract ad 2} again. 

Suppose { (\ref{cas des �clatements n1})\/} and let us show { (\ref{cas des �clatements n2})}. If $| r_{\a} r_{\b} | \neq \infty$, then necessarily $| r_{\a_k} r_{\b_k} |  \neq \infty$ for $k = 1,\, 2$. Moreover, $\G$ is then spherical (by proposition \ref{I spherique ssi tildeI sph�rique}), hence so is $\G_k$ for $k = 1,\, 2$. Lemma \ref{caract ad 2} gives us that $\gras r_{I} = \prod_{| r_{\a} r_{\b} |}(\gras r_\a,\gras r_\b) = \prod_{| r_{\a} r_{\b} |}(\gras r_\b,\gras r_\a)$. Let us denote by $I_k$ the vertex set of $\G_k$ for $k = 1,\, 2$. Since we have $\gras r_{I} = \gras r_{I_1}\gras r_{I_2}$, $\prod_{| r_{\a} r_{\b} |}(\gras r_\a,\gras r_\b) = \prod_{| r_{\a} r_{\b} |}(\gras r_{\a_1},\gras r_{\b_1})\prod_{| r_{\a} r_{\b} |}(\gras r_{\a_2},\gras r_{\b_2})$, and similarly if we exchange the roles of $\a$ and $\b$, and the roles of $a_k$ and $\b_k$, we conclude, by identifying the terms in $B^+_{\G_1}$ and $B^+_{\G_2}$, that $\gras r_{I_k} = \prod_{| r_{\a} r_{\b} |}(\gras r_{\a_k},\gras r_{\b_k}) = \prod_{| r_{\a} r_{\b} |}(\gras r_{\b_k},\gras r_{\a_k})$, for $k = 1,\, 2$. If $| r_{\a} r_{\b} | = \infty$, then lemma \ref{caract ad 2} shows us that the element $\prod_{m}(\gras r_\a,\gras r_\b) = \prod_{m}(\gras r_{\a_1},\gras r_{\b_1})\prod_{m}(\gras r_{\a_2},\gras r_{\b_2})$ is simple for all $m \in \Nset$, and similarly if we exchange the roles of $\a$ and $\b$, and the roles of $a_k$ and $\b_k$. We then have that $\prod_{m}(\gras r_{\a_k},\gras r_{\b_k})$ and $\prod_{m}(\gras r_{\b_k},\gras r_{\a_k})$ are simple for all $m \in \Nset$ (and $k = 1,\, 2$). In both cases, lemma \ref{LCM-donn�e et ordre} shows that $| r_\a r_\b | = | r_{\a_k} r_{\b_k} |$, for $k = 1,2$, and we conclude thanks to lemma \ref{caract ad 2}.\qed
\medskip

\subsection{Admissible $2$-partitions of irreducible spherical Coxeter graphs.}\label{parit� dans An}\mbox{}\medskip

Let $\G = (m_{i,j})_{i,j\in I}$ be a spherical Coxeter matrix and let $\tilde I = \{\a,\b\}$ be a $2$-partition of $\G$. Let us denote by $\G_1,\ldots,\, \G_n$ the connected components of $\G$, and by $I_k$ the vertex set of $\G_k$ for $1 \leqslant k \leqslant n$.

\begin{itemize}
\item If there exists $1 \leqslant k \leqslant n$ such that $I_k$ is included  in $\a$ or in $\b$, then $\tilde I$ is admissible if and only if $\G = \G_\a \times \G_\b$, in which case  $|r_\a r_\b| = 2$ (by proposition \ref{2.5.15}).
\item If not, then ($|r_\a r_\b| \not= 2$ and) $\a$ and $\b$ meet every connected component of $\G$, and we are in the situation of proposition \ref{cas des �clatements n}, with $\a_k = \a \cap I_k$ and $\b_k = \b \cap I_k$ for $1 \leqslant k \leqslant n$. So we get that $\tilde I = \{\a,\b\}$ is admissible if and only if $\{\a_k,\b_k\}$ is an admissible $2$-partition of $\G_{k}$ for $1 \leqslant k \leqslant n$, and $|r_\a r_\b| = |r_{\a_1} r_{\b_1}| = |r_{\a_2} r_{\b_2}| = \cdots = |r_{\a_n} r_{\b_n}|$. 
\end{itemize}

Hence we are left with the classification of admissible $2$-partitions of irreducible spherical Coxeter graphs and their corresponding coefficient $|r_\a r_\b|$. The first result in this direction in the following proposition. Since irreducible spherical Coxeter graphs are finite trees (hence bipartite), each of them has a unique {\em bipartite\/} partition, which is a $2$-partition except for the type $A_1$. 

\begin{prop}[{\cite[Lem.~2.5.13]{Mu2}}]\label{ad des partitions bipartites}The bipartite partition $\{\a, \b\}$ of an irreducible spherical Coxeter graph (distinct from $A_1$) is admissible, and the coefficient $|r_\a r_\b|$ is the {\em Coxeter number\/} of the graph.
\end{prop}

\proof  The result \cite[Lem.~5.8]{BS} gives our characterization of lemma \ref{caract ad 2}. \qed\medskip

\begin{rem}
 These considerations justify all cases of examples { \ref{sym spheriques}\/} and { \ref{H3 dans D6}\/} above except the ones concerning the non-trivial automorphisms of $A_{2n}$, $F_4$ and $I_2(m)$ (this last one being obvious), and reduce the justifications for $A_{2n}$ to the $A_4$ case. These last two cases (non-trivial automorphisms of $A_4$ and $F_4$) can be dealt with by direct computations.
\end{rem}

Let us now investigate the different situations case-by-case. 

\subsubsection{Admissible $2$-partitions of $A_n$, $B_n$, $D_n$.}\label{$2$-partitions admissibles des types An ...}\mbox{}\medskip

The admissible $2$-partitions of Coxeter graphs of type $A_n$ $(n\geqslant 2)$, $B_n$ $(n\geqslant 2)$ and $D_n$ $(n\geqslant 4)$ have been classified by M\"uhlherr in \cite[section~2.5]{Mu2}. In those cases, the only admissible $2$-partitions are the bipartite ones and, for every $n\geqslant 2$, the following $2$-partition of $A_{2n}$ (where the vertices are numbered in the natural order)~:  

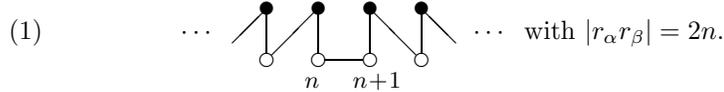
\begin{center}\label{$2$-partition exceptionnelle An}
\begin{picture}(60,45)(0,-15)\setlength{\unitlength}{1.3pt}
\put(0,0){\circle{4}}
\put(15,0){\circle{4}}\put(11,-8){\small{$n$}}
\put(30,0){\circle{4}}\put(25,-8){\small{$n\! +\! 1$}}
\put(45,0){\circle{4}}

\put(0,15){\circle*{4}}
\put(15,15){\circle*{4}}
\put(30,15){\circle*{4}}
\put(45,15){\circle*{4}}

\put(0,15){\line(-1,-1){10}}
\put(0,2){\line(0,1){14}}
\put(1.5,1.5){\line(1,1){14}}
\put(15,2){\line(0,1){14}}
\put(17,0){\line(1,0){11}}
\put(30,2){\line(0,1){14}}
\put(43.5,1.5){\line(-1,1){14}}
\put(45,2){\line(0,1){14}}
\put(45,15){\line(1,-1){10}}

\put(-25,6){$\cdots$}
\put(60,6){$\cdots$}
\put(75,6){with $|r_\a r_\b| = 2n$.}
\put(-75,6){(1)}
\end{picture}
\end{center}

The admissibility of this $2$-partition is a consequence of \cite[Lem.~2.5.6]{Mu2} (\cf{} proposition \ref{[M�2, 2.5.6]} above) applied to the admissible partition of $A_{2n}$ induced by its non-trivial automorphism and the bipartite partition of $B_n$. 

M\"uhlherr first established the classification for the $A_n$ case by explicit computations in the symmetric group. He inferred from this the classification for the $B_n$ case using \cite[Lem.~2.5.6]{Mu2}, which shows that every admissible $2$-partition of $B_n$ ``lifts'' to an admissible $2$-partition of $A_{2n}$ (or $A_{2n-1}$). In the same vein, since the automorphism of $D_n$ that permutes the vertices $n-1$ and $n$ (for the standard numbering of \cite[Planche~IV]{B}) gives an admissible partition of type $B_{n-1}$, and since \cite[Lem.~2.5.15]{Mu2} (\cf{} proposition \ref{2.5.15} above) shows that for every admissible $2$-partition of $D_n$, the vertices $n-1$ and $n$ must be in the same part of the partition, we get by \cite[Lem.~2.5.6]{Mu2} that every admissible $2$-partition of $D_n$ induces an admissible $2$-partition of $B_{n-1}$, whence the classification for the $D_n$ case. 

\subsubsection{Admissible $2$-partitions of $E_6$, $E_7$ and $E_8$.}\label{$2$-partitions admissibles des types En}\mbox{}\medskip

M\"uhlherr showed in \cite[Lem.~2.5.14]{Mu2} that the following $2$-partition of $E_6$ is admissible~: this is a consequence of \cite[Lem.~2.5.6]{Mu2} (\cf{} proposition \ref{[M�2, 2.5.6]} above) applied to the admissible partitions of $E_6$ and $F_4$ induced by their non-trivial automorphism. 

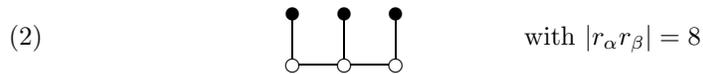
\begin{center}\label{$2$-partition exceptionnelle En}
\begin{picture}(60,40)(0,-10)\setlength{\unitlength}{1.3pt}
\put(7.5,0){\circle{4}}
\put(22.5,0){\circle{4}}
\put(37.5,0){\circle{4}}

\put(7.5,15){\circle*{4}}
\put(22.5,15){\circle*{4}}
\put(37.5,15){\circle*{4}}

\put(7.5,2){\line(0,1){14}}
\put(9.5,0){\line(1,0){11}}
\put(22.5,2){\line(0,1){14}}
\put(24.5,0){\line(1,0){11}}
\put(37.5,2){\line(0,1){14}}

\put(75,6){with $|r_\a r_\b| = 8$}
\put(-75,6){(2)}\label{part ad E6}
\end{picture}
\end{center}

\begin{prop}\label{classification En}The only admissible $2$-partitions of the Coxeter graphs $E_n$ $(n = 6,7,8)$ are the bipartite ones and the $2$-partition {\em (2)\/} above.
\end{prop}
\proof Let $\G$ be a Coxeter graph of type $E_6$, $E_7$ or $E_8$ and let $\{\a,\b\}$ be an admissible $2$-partition of $\G$. Since $\G$ is connected, $\{\a,\b\}$ does not satisfy the condition of proposition \ref{2.5.15} above. Hence, apart from the bipartite partitions and the $2$-partition (2) above, there are fifteen other possibilities~: 
\begin{itemize}
\item[-] one for $E_6$~: 
\end{itemize}
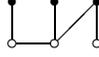
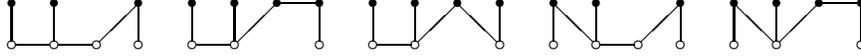
\begin{center}
\begin{picture}(45,15)(-7.5,0)\setlength{\unitlength}{0.8pt}
\put(0,0){\circle{4}}
\put(20,0){\circle{4}}
\put(40,0){\circle{4}}

\put(0,20){\circle*{4}}
\put(20,20){\circle*{4}}
\put(40,20){\circle*{4}}

\put(2,0){\line(1,0){16}}

\put(0,2){\line(0,1){17}}
\put(20,2){\line(0,1){17}}
\put(40,2){\line(0,1){17}}

\put(21.5,1.5){\line(1,1){18}}
\end{picture}
\end{center}
\begin{itemize}
\item[-] five for $E_7$~: 
\end{itemize}
\begin{center}
\begin{picture}(65,15)(-7.5,0)\setlength{\unitlength}{0.8pt}
\put(0,0){\circle{4}}
\put(20,0){\circle{4}}
\put(40,0){\circle{4}}
\put(60,0){\circle{4}}

\put(0,20){\circle*{4}}
\put(20,20){\circle*{4}}
\put(60,20){\circle*{4}}

\put(2,0){\line(1,0){16}}

\put(0,2){\line(0,1){17}}
\put(20,2){\line(0,1){17}}
\put(22,0){\line(1,0){16}}
\put(60,2){\line(0,1){17}}

\put(41.5,1.5){\line(1,1){18}}
\end{picture}
\begin{picture}(65,15)(-7.5,0)\setlength{\unitlength}{0.8pt}
\put(0,0){\circle{4}}
\put(20,0){\circle{4}}
\put(60,0){\circle{4}}

\put(0,20){\circle*{4}}
\put(20,20){\circle*{4}}
\put(40,20){\circle*{4}}
\put(60,20){\circle*{4}}

\put(2,0){\line(1,0){16}}

\put(0,2){\line(0,1){17}}
\put(20,2){\line(0,1){17}}
\put(42,20){\line(1,0){16}}
\put(60,2){\line(0,1){17}}

\put(21.5,1.5){\line(1,1){18}}
\end{picture}
\begin{picture}(65,15)(-7.5,0)\setlength{\unitlength}{0.8pt}
\put(0,0){\circle{4}}
\put(20,0){\circle{4}}
\put(60,0){\circle{4}}

\put(0,20){\circle*{4}}
\put(20,20){\circle*{4}}
\put(40,20){\circle*{4}}
\put(60,20){\circle*{4}}

\put(2,0){\line(1,0){16}}

\put(0,2){\line(0,1){17}}
\put(20,2){\line(0,1){17}}
\put(58.5,1.5){\line(-1,1){18}}
\put(60,2){\line(0,1){17}}

\put(21.5,1.5){\line(1,1){18}}
\end{picture}
\begin{picture}(65,15)(-7.5,0)\setlength{\unitlength}{0.8pt}
\put(0,0){\circle{4}}
\put(20,0){\circle{4}}
\put(40,0){\circle{4}}
\put(60,0){\circle{4}}

\put(0,20){\circle*{4}}
\put(20,20){\circle*{4}}
\put(60,20){\circle*{4}}

\put(0,2){\line(0,1){16}}

\put(18.5,1.5){\line(-1,1){18}}
\put(20,2){\line(0,1){17}}
\put(22,0){\line(1,0){16}}
\put(60,2){\line(0,1){17}}

\put(41.5,1.5){\line(1,1){18}}
\end{picture}
\begin{picture}(65,15)(-7.5,0)\setlength{\unitlength}{0.8pt}
\put(0,0){\circle{4}}
\put(20,0){\circle{4}}
\put(60,0){\circle{4}}

\put(0,20){\circle*{4}}
\put(20,20){\circle*{4}}
\put(40,20){\circle*{4}}
\put(60,20){\circle*{4}}

\put(18.5,1.5){\line(-1,1){18}}

\put(0,2){\line(0,1){17}}
\put(20,2){\line(0,1){17}}
\put(42,20){\line(1,0){16}}
\put(60,2){\line(0,1){17}}

\put(21.5,1.5){\line(1,1){18}}
\end{picture}
\end{center}
\begin{itemize}
\item[-] and nine for $E_8$~: 
\end{itemize}
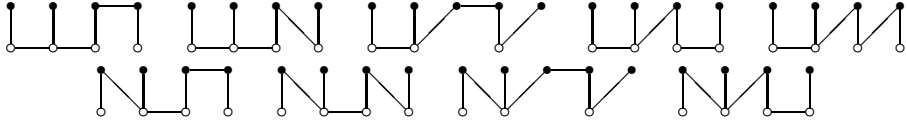
\begin{center}
\begin{picture}(65,15)(0,0)\setlength{\unitlength}{0.8pt}
\put(0,0){\circle{4}}
\put(20,0){\circle{4}}
\put(40,0){\circle{4}}
\put(60,0){\circle{4}}

\put(0,20){\circle*{4}}
\put(20,20){\circle*{4}}
\put(40,20){\circle*{4}}
\put(60,20){\circle*{4}}

\put(2,0){\line(1,0){16}}
\put(42,20){\line(1,0){16}}

\put(0,2){\line(0,1){16}}
\put(20,2){\line(0,1){16}}
\put(40,2){\line(0,1){16}}
\put(60,2){\line(0,1){16}}

\put(22,0){\line(1,0){16}}
\end{picture}
\begin{picture}(65,15)(0,0)\setlength{\unitlength}{0.8pt}
\put(0,0){\circle{4}}
\put(20,0){\circle{4}}
\put(40,0){\circle{4}}
\put(60,0){\circle{4}}

\put(0,20){\circle*{4}}
\put(20,20){\circle*{4}}
\put(40,20){\circle*{4}}
\put(60,20){\circle*{4}}

\put(2,0){\line(1,0){16}}
\put(58.5,1.5){\line(-1,1){18}}

\put(0,2){\line(0,1){16}}
\put(20,2){\line(0,1){16}}
\put(40,2){\line(0,1){16}}
\put(60,2){\line(0,1){16}}

\put(22,0){\line(1,0){16}}
\end{picture}
\begin{picture}(80,15)(0,0)\setlength{\unitlength}{0.8pt}
\put(0,0){\circle{4}}
\put(20,0){\circle{4}}
\put(80,20){\circle*{4}}
\put(60,0){\circle{4}}

\put(0,20){\circle*{4}}
\put(20,20){\circle*{4}}
\put(40,20){\circle*{4}}
\put(60,20){\circle*{4}}

\put(2,0){\line(1,0){16}}
\put(42,20){\line(1,0){16}}

\put(0,2){\line(0,1){16}}
\put(20,2){\line(0,1){16}}
\put(60,2){\line(0,1){16}}
\put(61.5,1.5){\line(1,1){18}}
\put(21.5,1.5){\line(1,1){18}}
\end{picture}
\begin{picture}(65,15)(0,0)\setlength{\unitlength}{0.8pt}
\put(0,0){\circle{4}}
\put(20,0){\circle{4}}
\put(40,0){\circle{4}}
\put(60,0){\circle{4}}

\put(0,20){\circle*{4}}
\put(20,20){\circle*{4}}
\put(40,20){\circle*{4}}
\put(60,20){\circle*{4}}

\put(2,0){\line(1,0){16}}
\put(42,0){\line(1,0){16}}

\put(0,2){\line(0,1){16}}
\put(20,2){\line(0,1){16}}
\put(40,2){\line(0,1){16}}
\put(60,2){\line(0,1){16}}

\put(21.5,1.5){\line(1,1){18}}
\end{picture}
\begin{picture}(65,15)(0,0)\setlength{\unitlength}{0.8pt}
\put(0,0){\circle{4}}
\put(20,0){\circle{4}}
\put(40,0){\circle{4}}
\put(60,0){\circle{4}}

\put(0,20){\circle*{4}}
\put(20,20){\circle*{4}}
\put(40,20){\circle*{4}}
\put(60,20){\circle*{4}}

\put(2,0){\line(1,0){16}}

\put(0,2){\line(0,1){16}}
\put(20,2){\line(0,1){16}}
\put(40,2){\line(0,1){16}}
\put(60,2){\line(0,1){16}}
\put(41.5,1.5){\line(1,1){18}}
\put(21.5,1.5){\line(1,1){18}}
\end{picture}
\end{center}
\begin{center}
\begin{picture}(65,15)(0,0)\setlength{\unitlength}{0.8pt}
\put(0,0){\circle{4}}
\put(20,0){\circle{4}}
\put(40,0){\circle{4}}
\put(60,0){\circle{4}}

\put(0,20){\circle*{4}}
\put(20,20){\circle*{4}}
\put(40,20){\circle*{4}}
\put(60,20){\circle*{4}}

\put(18.5,1.5){\line(-1,1){18}}
\put(42,20){\line(1,0){16}}

\put(0,2){\line(0,1){16}}
\put(20,2){\line(0,1){16}}
\put(40,2){\line(0,1){16}}
\put(60,2){\line(0,1){16}}

\put(22,0){\line(1,0){16}}
\end{picture}
\begin{picture}(65,15)(0,0)\setlength{\unitlength}{0.8pt}
\put(0,0){\circle{4}}
\put(20,0){\circle{4}}
\put(40,0){\circle{4}}
\put(60,0){\circle{4}}

\put(0,20){\circle*{4}}
\put(20,20){\circle*{4}}
\put(40,20){\circle*{4}}
\put(60,20){\circle*{4}}

\put(18.5,1.5){\line(-1,1){18}}
\put(58.5,1.5){\line(-1,1){18}}

\put(0,2){\line(0,1){16}}
\put(20,2){\line(0,1){16}}
\put(40,2){\line(0,1){16}}
\put(60,2){\line(0,1){16}}

\put(22,0){\line(1,0){16}}
\end{picture}
\begin{picture}(80,15)(0,0)\setlength{\unitlength}{0.8pt}
\put(0,0){\circle{4}}
\put(20,0){\circle{4}}
\put(80,20){\circle*{4}}
\put(60,0){\circle{4}}

\put(0,20){\circle*{4}}
\put(20,20){\circle*{4}}
\put(40,20){\circle*{4}}
\put(60,20){\circle*{4}}

\put(18.5,1.5){\line(-1,1){18}}
\put(42,20){\line(1,0){16}}

\put(0,2){\line(0,1){16}}
\put(20,2){\line(0,1){16}}
\put(60,2){\line(0,1){16}}
\put(61.5,1.5){\line(1,1){18}}
\put(21.5,1.5){\line(1,1){18}}
\end{picture}
\begin{picture}(65,15)(0,0)\setlength{\unitlength}{0.8pt}
\put(0,0){\circle{4}}
\put(20,0){\circle{4}}
\put(40,0){\circle{4}}
\put(60,0){\circle{4}}

\put(0,20){\circle*{4}}
\put(20,20){\circle*{4}}
\put(40,20){\circle*{4}}
\put(60,20){\circle*{4}}

\put(18.5,1.5){\line(-1,1){18}}
\put(42,0){\line(1,0){16}}

\put(0,2){\line(0,1){16}}
\put(20,2){\line(0,1){16}}
\put(40,2){\line(0,1){16}}
\put(60,2){\line(0,1){16}}

\put(21.5,1.5){\line(1,1){18}}
\end{picture}
\end{center}

By lemma \ref{caract ad}, there exist $n \in \Nset$ such that $\prod_{n}(r_\a, r_\b) = \prod_{n}(r_\b, r_\a) = r_I$ and $\sum_{n}(\l(r_\a), \l(r_\b)) = \sum_{n}(\l(r_\b), \l(r_\a)) = \l(r_I)$. Since we have $\l(r_I) = 36$ (resp. $63$, resp. $120$) if $\G = E_6$ (resp. $E_7$, resp. $E_8$), \cf{} \cite[Planches~V-VII]{B}, the consideration  on lengths eliminates the last candidate for $E_6$ and leaves only one candidate for $E_7$ (the second one, with $n = 14$) and four for $E_8$ (the first one with $n = 20$, and the third, fourth and sixth ones with $n = 24$). We then verify, if needed with the help of a computation software like GAP or Maple, that the equality $\prod_{n}(r_\a, r_\b) = \prod_{n}(r_\b, r_\a) = r_I$ occurs in none of the five remaining cases, hence those $2$-partitions are not admissible. \qed\medskip

\subsubsection{Admissible $2$-partitions of $F_4$, $H_3$, $H_4$ (and $I_2(m)$, $m\geqslant 3$).}\label{$2$-partitions admissibles des types F4 ...}\mbox{}\medskip

The orbits of $F_4$ under the action of its non-trivial automorphism form the following admissible $2$-partition~: 

\begin{center}\label{$2$-partition exceptionnelle Fn}
\begin{picture}(60,40)(0,-10)\setlength{\unitlength}{1.3pt}
\put(15,0){\circle{4}}
\put(30,0){\circle{4}}

\put(15,15){\circle*{4}}
\put(30,15){\circle*{4}}

\put(15,2){\line(0,1){14}}
\put(17,0){\line(1,0){11}}\put(21,1.5){\small{4}}
\put(30,2){\line(0,1){14}}

\put(75,6){with $|r_\a r_\b| = 8$}
\put(-75,6){(3)}
\end{picture}
\end{center}

\begin{prop}\label{classification F4 ...}The only admissible $2$-partitions of the Coxeter graphs $F_4$, $H_3$, $H_4$ and $I_2(m)$, $m\geqslant 3$, are the bipartite ones and the $2$-partition {\em (3)\/} above.
\end{prop}
\proof There is nothing to prove for the dihedral graphs. So assume that $\G$ is a Coxeter graph of type $F_4$, $H_3$ or $H_4$, and let $\{\a,\b\}$ be an admissible $2$-partition of $\G$. Since $\G$ is connected, $\{\a,\b\}$ does not satisfy the condition of proposition \ref{2.5.15} above and hence is either a bipartite partition, or the $2$-partition (3) above, or possibly the following $2$-partition of $H_4$~: 
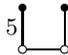
\begin{center}
\begin{picture}(20,15)(0,0)\setlength{\unitlength}{0.8pt}
\put(0,0){\circle{4}}
\put(20,0){\circle{4}}

\put(0,20){\circle*{4}}
\put(20,20){\circle*{4}}

\put(0,2){\line(0,1){18}}\put(-8,7){\small{5}}
\put(2,0){\line(1,0){16}}
\put(20,2){\line(0,1){18}}
\end{picture}
\end{center}

To show that this last $2$-partition is non-admissible, one can follow the same lines as in the proof of proposition \ref{classification En}. Otherwise, note that $H_4$ is the type of an admissible partition of $E_8$ (\cf{} examples \ref{H3 dans D6} or \ref{eclatement H3 H4}) so, thanks to proposition \ref{[M�2, 2.5.6]}, the admissibility of the above $2$-partition of $H_4$ is equivalent to the admissibility of a certain $2$-partition of $E_8$ (not the bipartite one), which has been shown to be non-admissible in proposition \ref{classification En}. \qed\medskip

\subsection{Foldings.}\label{Pliages.}\mbox{}\medskip

Let $\G = (m_{i,j})_{i,j\in I}$ and $\G' = (m'_{i',j'})_{i',j'\in I'}$ be two Coxeter matrices with no infinite coefficient. Crisp defined in \cite[Def.~4.1]{C} the notion of a {\em folding\/} of $\G'$ onto $\G$, in order to give examples of LCM-homomorphisms and to begin their classification. With our terminology, a folding of $\G'$ onto $\G$ is a surjective map $f : I' \twoheadrightarrow I$ that satisfy a list of conditions made for the partition $\{f^{-1}(\{i\}) \mid i\in I\}$ of $I'$ to be an LCM-partition of type (isomorphic to) $\G$ \cite[Prop.~4.2]{C}. Crisp concluded \cite[section~4]{C} by asking essentially if every LCM-partition is obtained from a folding. The classification we have just established shows that the answer is no, with the definition \cite[Def.~4.1]{C} for a folding, and indicates how to complete the list of cases of \cite[Def.~4.1]{C} to turn the answer to yes. 

In definition \ref{folding} below, we propose a generalisation of the notion of foldings that fit to our new point of view, and in proposition \ref{r�ponse Crisp}, we rephrase in the manner of \cite[Def.~4.1]{C} the classification established above.

\begin{defn}[foldings]\label{folding}Let $\G = (m_{i,j})_{i,j\in I}$ and $\G' = (m'_{i',j'})_{i',j'\in I'}$ be two Coxeter matrices. A {\em folding\/} of $\G'$ onto $\G$ is a map $f : I' \to I$ such that the set $\{f^{-1}(\{i\}) \mid i \in I\}$ is an admissible partition of $\G'$, of type (isomorphic to) $\G$.
\end{defn}

\begin{nota}Let $\G$ be any Coxeter graph. For $n \in \Nset_{\geqslant 1}$, we denote by $n\G$ the disjoint union of $n$ copies of $\G$. 
\end{nota}

\begin{prop}\label{r�ponse Crisp}Let $\G = (m_{i,j})_{i,j\in I}$ and $\G' = (m'_{i',j'})_{i',j'\in I'}$ be two Coxeter matrices and $f : I' \to I$ be any map from $I'$ to $I$. Assume that $\G$ has no infinite coefficient. Then $f$ is a folding from $\G'$ onto $\G$ if and only if $f$ satisfies the following conditions for every $i, j \in I$~: 
\begin{enumerate}
\item \label{r�ponse Crisp1} the subset $f^{-1}(\{i\})$ of $I'$ is non-empty and spherical,
\item \label{r�ponse Crisp2} if $m_{i,j} = 2$, then there is no edge between a vertex of $f^{-1}(\{i\})$ and a vertex of $f^{-1}(\{j\})$, \ie{} $\G'_{f^{-1}(\{i,j\})} = \G'_{f^{-1}(\{i\})} \times \G'_{f^{-1}(\{j\})}$,
\item \label{r�ponse Crisp3} if $m_{i,j} \geqslant 3$, then one of the following occurs~: 
\begin{enumerate}
\item[\bf(A)] $\G'_{f^{-1}(\{i,j\})} = nI_2(m_{i,j})$ for some $n\in \Nset_{\geqslant 1}$, and each connected component  of $\G'_{f^{-1}(\{i,j\})}$ (of type $I_2(m_{i,j})$) meets $f^{-1}(\{i\})$ and $f^{-1}(\{j\})$,
\item[\bf(B)] $\G'_{f^{-1}(\{i,j\})}$ is an irreducible and spherical Coxeter graph with Coxeter number $m_{i,j}$, and the $2$-partition $\{f^{-1}(\{i\}), f^{-1}(\{j\})\}$ of $f^{-1}(\{i,j\})$ is the bipartite partition of $\G'_{f^{-1}(\{i,j\})}$,
\item[\bf(C1)] $m_{i,j} = 2n$ for some $n \in \Nset_{\geqslant 2}$, $\G'_{f^{-1}(\{i,j\})} = A_{2n}$, and the $2$-partition $\{f^{-1}(\{i\}), f^{-1}(\{j\})\}$ of $f^{-1}(\{i,j\})$ is the admissible $2$-partition {\em (1)\/} of subsection {\em \ref{$2$-partitions admissibles des types An ...}},
\item[\bf(C2)] $m_{i,j} = 8$, $\G'_{f^{-1}(\{i,j\})} = E_6$, and the $2$-partition $\{f^{-1}(\{i\}), f^{-1}(\{j\})\}$ of $f^{-1}(\{i,j\})$ is the admissible $2$-partition {\em (2)\/} of subsection {\em \ref{$2$-partitions admissibles des types En}},
\item[\bf(C3)] $m_{i,j} = 8$, $\G'_{f^{-1}(\{i,j\})} = F_4$, and the $2$-partition $\{f^{-1}(\{i\}), f^{-1}(\{j\})\}$ of $f^{-1}(\{i,j\})$ is the admissible $2$-partition {\em (3)\/} of subsection {\em \ref{$2$-partitions admissibles des types F4 ...}},
\item[\bf(D)] the map $f^{-1}(\{i,j\}) \to \{i,j\}$ induced by $f$ is a composition $h\circ g$, where $g$ is a folding from $\G'_{f^{-1}(\{i,j\})}$ onto $nI_2(m_{i,j})$ ($n\in \Nset_{\geqslant 2}$) defined only with cases {\em (B)\/} to {\em (C3)\/} and $h$ is a folding from $nI_2(m_{i,j})$ onto $\G_{\{i,j\}} = I_2(m_{i,j})$ of case {\em (A)}.
\end{enumerate}
\end{enumerate}
\end{prop}
\proof This is a reformulation of the classification obtained above.\qed\medskip

\begin{rem}We have added to the list of { \cite[Def.~4.1]{C}\/} the cases { (C1)\/} for $n > 2$, { (C2)} and { (C3)\/}. Note that  { \cite[Def.~1.11]{G2}\/} already includes case { (C3)}.
\end{rem}

\begin{rem}The cases { (A)\/} to  { (D)\/} imply, for a non-isolated vertex $i$ of $\G$, that $\G'_{f^{-1}(\{i\})}$ is non-empty and spherical, hence our condition (\ref{r�ponse Crisp1}) can be relaxed to the weaker condition (implicit in { \cite[Def.~4.1]{C}\/} and { \cite[Def.~1.11]{G2}})~:
\begin{enumerate}
\item[(\ref{r�ponse Crisp1}')] if $i$ is isolated in $\G$, then $f^{-1}(\{i\})$ is non-empty and spherical.
\end{enumerate}
\end{rem}

\end{document}